\documentclass[12pt,vatola]{article}

\usepackage{graphicx}

\def\c{\centerline}

\def\no{\noindent}

\begin{document}

\c{\Large\bf Smarandache Multi-Space Theory(I) }\vskip 5mm

\hskip 70mm {\it -Algebraic multi-spaces }\vskip 10mm

\c{Linfan Mao}\vskip 3mm

\c{\small Academy of Mathematics and System Sciences}

\c{\small Chinese Academy of Sciences, Beijing 100080}

\c{\small maolinfan@163.com}

\vskip 10mm

\begin{minipage}{130mm}

\no{\bf Abstract.} {\small A Smarandache multi-space is a union of
$n$ different spaces equipped with some different structures for
an integer $n\geq 2$, which can be both used for discrete or
connected spaces, particularly for geometries and spacetimes in
theoretical physics. This monograph concentrates on characterizing
various multi-spaces including three parts altogether. The first
part is on {\it algebraic multi-spaces with structures}, such as
those of multi-groups, multi-rings, multi-vector spaces,
multi-metric spaces, multi-operation systems and multi-manifolds,
also multi-voltage graphs, multi-embedding of a graph in an
$n$-manifold,$\cdots$, etc.. The second discusses {\it Smarandache
geometries}, including those of map geometries, planar map
geometries and pseudo-plane geometries, in which the {\it Finsler
geometry}, particularly the {\it Riemann geometry} appears as a
special case of these Smarandache geometries. The third part of
this book considers the {\it applications of multi-spaces to
theoretical physics}, including the relativity theory, the
M-theory and the cosmology. Multi-space models for $p$-branes and
cosmos are constructed and some questions in cosmology are
clarified by multi-spaces. The first two parts are relative
independence for reading and in each part open problems are
included for further research of interested readers.}

\vskip 3mm \no{\bf Key words:} {\small  algebraic structure,
multi-space, multi-group, multi-ring, multi-vector space,
multi-metric space.}

 \vskip 3mm \no{{\bf
Classification:} AMS(2000) 03C05,05C15,51D20,51H20,51P05,83C05,
83E50}
\end{minipage}

\newpage

\no{\large\bf Contents}\vskip 8mm

\no $1.$\dotfill 3\vskip 3mm

\no $\S 1.1$ \ sets\dotfill 3\vskip 2mm

\no $1.1.1$ sets \dotfill 3

\no $1.1.2$ Partially ordered sets  \dotfill 5

\no $1.1.3$ Neutrosophic sets  \dotfill 7\vskip 2mm

\no $\S 1.2$ \ Algebraic Structures \dotfill 8\vskip 2mm

\no $1.2.1$ Groups \dotfill 9

\no $1.2.2$ Rings \dotfill 10

\no $1.2.3$ Vector spaces \dotfill 12\vskip 2mm

\no $\S 1.3$ \ Algebraic Multi-Spaces\dotfill 14\vskip 2mm

\no $1.3.1$ Algebraic multi-spaces\dotfill 14

\no $1.3.2$ Multi-Groups \dotfill 20

\no $1.3.3$ Multi-Rings \dotfill 25

\no $1.3.4$ Multi-Vector spaces \dotfill 31\vskip 2mm

\no $\S 1.4$ \ Multi-Metric Spaces\dotfill 35\vskip 2mm

\no $1.4.1$ Metric spaces\dotfill 35

\no $1.4.2$ Multi-Metric spaces\dotfill 36\vskip 2mm

\no $\S 1.5$ \ Remarks and Open Problems\dotfill 42\vskip 3mm

\newpage

\no{\large\bf $1.$ Algebraic multi-spaces}\vskip 25mm

The notion of multi-spaces was introduced by Smarandache in 1969,
see his article uploaded to {\it arXiv} $[86]$ under his idea of
hybrid mathematics: {\it combining different fields into a
unifying field}($[85]$), which is more closer to our real life
world. Today, this idea is widely accepted by the world of
sciences. For mathematics, a definite or an exact solution under a
given condition is not the only object for mathematician. New
creation power has emerged and new era for mathematics has come
now. Applying the Smarandache's notion, this chapter concentrates
on constructing various multi-spaces by algebraic structures, such
as those of groups, rings, fields, vector spaces, $\cdots$,etc.,
also by metric spaces, which are more useful for constructing
multi-voltage graphs, maps and map geometries in the following
chapters.

\vskip 8mm

\no{\bf \S $1.1$ \ Sets}

\vskip 5mm

\no{\bf $1.1.1.$ Sets}

\vskip 4mm

\no A {\it set} $\Xi$ is a collection of objects with some common
property $P$, denoted by

$$\Xi =\{x| x \ has \ property \ P\},$$

\no where, $x$ is said an element of the set $\Xi$, denoted by
$x\in \Xi$. For an element $y$ not possessing the property $P$,
i.e., not an element in the set $\Xi$, we denote it by $y\not\in
\Xi$.

The {\it cardinality} (or the number of elements if $\Xi$ is
finite ) of a set $\Xi$ is denoted by $|\Xi|$.

Some examples of sets are as follows.

$$A=\{1,2,3,4,5,6,7,8,9,10\};$$

$$B=\{ p | \ p \ is \ a \ prime \ number \};$$

$$C=\{(x,y)| x^2+y^2=1\};$$

$$D=\{ the \ cities \ in \ the \ world \}.$$

The sets $A$ and $D$ are finite with $|A|=10$ and $|D| < +\infty$,
but these sets $B$ and $C$ are infinite.

Two sets $\Xi_1$ and $\Xi_2$ are said to be {\it identical} if and
only if for $\forall x\in\Xi_1$, we have $x\in\Xi_2$ and for
$\forall x\in\Xi_2$, we also have $x\in\Xi_1$. For example, the
following two sets

$$E=\{1,2,-2\} \ {\rm and} \ F=\{ \ x \ | x^3-x^2-4x+4=0\}$$

\no are identical since we can solve the equation $x^3-x^2-4x+4=0$
and get the solutions $x=1,2$ or $-2$. Similarly, for the
cardinality of a set, we know the following result.

\vskip 4mm

\no{\bf Theorem $1.1.1$}([6]) \ {\it For sets $\Xi_1, \Xi_2$,
$|\Xi_1|= |\Xi_2|$ if and only if there is an $1-1$ mapping
between $\Xi_1$ \ and \ $\Xi_2$.}

\vskip 3mm

According to this theorem, we know that $|B|\not=|C|$ although
they are infinite. Since $B$ is countable, i.e., there is an $1-1$
mapping between $B$ and the natural number set $N=\{1,2,3,\cdots
,n,\cdots\}$, however $C$ is not.

Let $A_1,A_2$ be two sets. If for $\forall a\in A_1\Rightarrow
a\in A_2$, then $A_1$ is said to be a {\it subset} of $A_2$,
denoted by $A_1\subseteq A_2$. If a set has no elements, we say it
an empty set, denoted by $\emptyset$.

\vskip 4mm

\no{\bf Definition $1.1.1$} \ {\it For two sets $\Xi_1,\Xi_2$, two
operations ¡°$\bigcup$¡±and ¡°$\bigcap$¡± on $\Xi_1,\Xi_2$ are
defined as follows:}

$$\Xi_1\bigcup\Xi_2=\{x|x\in\Xi_1 \ or \ x\in\Xi_2\},$$

$$\Xi_1\bigcap\Xi_2=\{x|x\in\Xi_1 \ and \ x\in\Xi_2\}$$

\no{\it and $\Xi_1$ minus $\Xi_2$ is defined by}

$$\Xi_1\setminus\Xi_2=\{x|x\in\Xi_1 \ but \ x\not\in\Xi_2\}.$$

\vskip 3mm

For the sets $A$ and $E$, calculation shows that

$$A\bigcup E=\{1,2,-2,3,4,5,6,7,8,9,10\},$$

$$A\bigcap E=\{1,2\}$$

\no and

$$A\setminus E=\{3,4,5,6,7,8,9,10\},$$

$$E\setminus A=\{-2\}.$$

For a set $\Xi$ and $H\subseteq\Xi$, the set $\Xi\setminus H$ is
said the {\it complement} of $H$ in $\Xi$, denoted by
$\overline{H}(\Xi )$. We also abbreviate it to $\overline{H}$ if
each set considered in the situation is a subset of $\Xi =
\Omega$, i.e., the {\it universal set}.

These operations defined in Definition $1.1.1$ observe the
following laws.

\vskip 2mm

\no{\bf L1} \ Itempotent law. \ For $\forall S\subseteq \Omega$,

$$A\bigcup A= A\bigcap A=A.$$

\vskip 2mm

\no{\bf L2} \ Commutative law. \ For $\forall U,V\subseteq
\Omega$,

$$U\bigcup V=V\bigcup U; \ U\bigcap V=V\bigcap U.$$

\vskip 2mm

\no{\bf L3} \ Associative law. \  For $\forall U,V, W\subseteq
\Omega$,

$$U\bigcup (V\bigcup W)= (U\bigcup V)\bigcup W; \
U\bigcap (V\bigcap W)=(U\bigcap V)\bigcap W.$$

\vskip 2mm

\no{\bf L4} \ Absorption law. \ \ For $\forall U,V\subseteq
\Omega$,

$$U\bigcap (U\bigcup V)=U\bigcup (U\bigcap V)=U.$$

\vskip 2mm

\no{\bf L5} \ Distributive law. \  For $\forall U,V, W\subseteq
\Omega$,

$$U\bigcup (V\bigcap W)= (U\bigcup V)\bigcap (U\bigcup W); \
U\bigcap (V\bigcup W)=(U\bigcap V)\bigcup (U\bigcap W).$$

\vskip 2mm

\no{\bf L6} \ Universal bound law. \ For $\forall U\subseteq
\Omega,$

$$\emptyset\bigcap U=\emptyset,\emptyset\bigcup U=U; \ \Omega\bigcap U=U, \Omega\bigcup U=\Omega.$$

\vskip 2mm

\no{\bf L7} \ Unary complement law. \ For $\forall U\subseteq
\Omega,$

$$U\bigcap\overline{U}=\emptyset ; \ U\bigcup\overline{U}=\Omega.$$

\vskip 2mm

A set with two operations ¡°$\bigcap$¡± and ¡°$\bigcup$¡±
satisfying the laws $L1\sim L7$ is said to be a {\it Boolean
algebra}. Whence, we get the following result.

\vskip 4mm

\no{\bf Theorem $1.1.2$} \ {\it For any set $U$, all its subsets
form a Boolean algebra under the operations ¡°$\bigcap$¡± and
¡°$\bigcup$¡±.}

\vskip 5mm

\no{\bf $1.1.2$ \ Partially order sets}

\vskip 3mm

\no For a set $\Xi$, define its {\it Cartesian product} to be

$$\Xi\times\Xi =\{(x,y)| \forall x,y\in\Xi\}.$$

A subset $R\subseteq\Xi\times\Xi$ is called a {\it binary
relation} on $\Xi$. If $(x,y)\in R$, we write $xRy$. A {\it
partially order set} is a set $\Xi$ with a binary relation
¡°$\preceq$¡± such that the following laws hold.

\vskip 2mm

\no{\bf O1} Reflective law. \ For $x\in\Xi$, \ $xRx$.

\vskip 2mm

\no{\bf O2} Antisymmetry law. \ For $x,y\in\Xi$, $xRy$ and
$yRx\Rightarrow x=y.$

\vskip 2mm

\no{\bf O3} Transitive law. \ For $x,y,z\in\Xi$, $xRy$ and
$yRz\Rightarrow xRz.$

\vskip 2mm

A partially order set $\Xi$ with a binary relation \ ¡°$\preceq$¡±
\ is denoted by $(\Xi ,\preceq )$. Partially ordered sets with a
finite number of elements can be conveniently represented by a
diagram in such a way that each element in the set $\Xi$ is
represented by a point so placed on the plane that a point $a$ is
above another point $b$ if and only if $b\prec a$. This kind of
diagram is essentially a directed graph (see also Chapter $2$ in
this book). In fact, a directed graph is correspondent with a
partially set and vice versa. Examples for the partially order
sets are shown in Fig.$1.1$ where each diagram represents a finite
partially order set.

\vskip 2mm

\includegraphics[bb=-30 10 200 110]{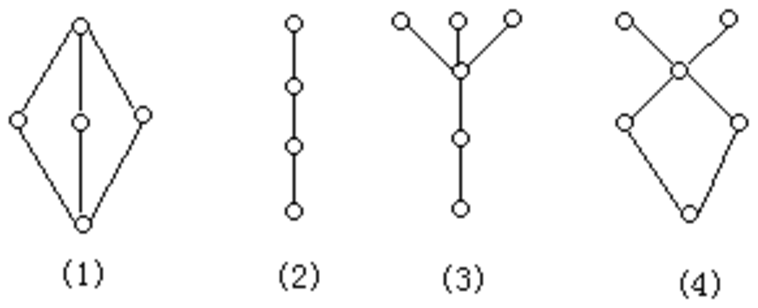}

\vskip 2mm

\c{\bf Fig.$1.1$}

\vskip 2mm

An element $a$ in a partially order set $(\Xi ,\preceq )$ is
called {\it maximal} (or {\it minimal}) if for $\forall x\in\Xi, \
a\preceq x\Rightarrow x=a$ (or $x\preceq a\Rightarrow x=a$). The
following result is obtained by the definition of partially order
sets and the induction principle.

\vskip 4mm

\no{\bf Theorem $1.1.3$} \ {\it Any finite non-empty partially
order set $(\Xi ,\preceq )$ has maximal and minimal elements.}

\vskip 3mm

A partially order set $(\Xi ,\preceq )$ is an {\it order set} if
for any $\forall x,y\in\Xi$, there must be $x\preceq y$ or
$y\preceq x$. It is obvious that any partially order set contains
an order subset, finding this fact in Fig.$1.1$.

An {\it equivalence relation} $R\subseteq\Xi\times\Xi$ on a set
$\Xi$ is defined by

\vskip 2mm

\no{\bf R1} Reflective law. \ For $x\in\Xi$, \ $xRx$.

\vskip 2mm

\no{\bf R2} Symmetry law. \ For $x,y\in\Xi$, $xRy \Rightarrow yRx
$

\vskip 2mm

\no{\bf R3} Transitive law. \ For $x,y,z\in\Xi$, $xRy$ and
$yRz\Rightarrow xRz.$

\vskip 2mm

For a set $\Xi$ with an equivalence relation $R$, we can classify
elements in $\Xi$ by $R$ as follows:

$$R(x)=\{ y | \ y\in\Xi \ {\rm and} \ yRx \ \}.$$

Then, we get the following useful result for the combinatorial
enumeration.

\vskip 4mm

\no{\bf Theorem $1.1.4$} \ {\it For a finite set $\Xi$ with an
equivalence $R$, $\forall x,y\in\Xi$, if there is an bijection
$\varsigma$ between $R(x)$ and $R(y)$, then the number of
equivalence classes under $R$ is}

$$\frac{|\Xi |}{|R(x)|},$$

\no{\it where $x$ is a chosen element in $\Xi$.}

\vskip 3mm

{\it Proof} \ Notice that there is an bijection $\varsigma$
between $R(x)$ and $R(y)$ for $\forall x,y\in\Xi$. Whence,
$|R(x)|=|R(y)|$. By definition, for $\forall x,y\in\Xi$,
$R(x)\bigcap R(y)=\emptyset$ or $R(x)=R(y)$. Let $T$ be a
representation set of equivalence classes, i.e., choice one
element in each class. Then we get that

\begin{eqnarray*}
|\Xi | &=& \sum\limits_{x\in T}|R(x)|\\
&=& |T||R(x)|.
\end{eqnarray*}

\no Whence, we know that

$$|T| \ = \ \frac{|\Xi |}{|R(x)|}.\quad\quad\quad\quad \natural$$

\vskip 5mm

\no{\bf $1.1.3$ \ Neutrosophic set}

\vskip 3mm

\no Let $[0,1]$ be a closed interval. For three subsets $T, I,
F\subseteq [0,1]$ and $S\subseteq \Omega$, define a relation of an
element $x\in\Omega$ with the subset $S$ to be $x(T,I,F)$, i,e.,
the {\it confidence set} for $x\in S$ is $T$, the {\it indefinite
set} is $I$ and {\it fail set} is $F$. A set $S$ with three
subsets $T,I,F$ is said to be a {\it neutrosophic set} ([85]). We
clarify the conception of neutrosophic sets by abstract set theory
as follows.

Let $\Xi$ be a set and $A_1,A_2,\cdots , A_k\subseteq\Xi$. Define
$3k$ functions $f_1^x,f_2^x,\cdots ,f_k^x$ by $f_i^x:
A_i\rightarrow [0,1], \  1\leq i\leq k$, where $x=T, I, F$. Denote
by $(A_i;f_i^T,f_i^I,f_i^F)$ the subset $A_i$ with three functions
$f_i^T, f_i^I, f_i^F$, $1\leq i\leq k$. Then

$$\bigcup\limits_{i=1}^k(A_i;f_i^T,f_i^I,f_i^F)$$

\no is a union of neutrosophic sets. Some extremal cases for this
union is in the following, which convince us that neutrosophic
sets are a generalization of classical sets.

\vskip 2mm

\no{\bf Case $1$} \ \ $f_i^T=1, \ f_i^I=f_i^F=0$ for $i=1,2,\cdots
,k$.

\vskip 2mm

In this case,

$$\bigcup\limits_{i=1}^k(A_i;f_i^T,f_i^I,f_i^F)=\bigcup\limits_{i=1}^kA_i.$$

\vskip 2mm

\no{\bf Case $2$} \ \ $f_i^T=f_i^I=0, \ f_i^F=1$ for $i=1,2,\cdots
, k$.

\vskip 2mm

In this case,

$$\bigcup\limits_{i=1}^k(A_i;f_i^T,f_i^I,f_i^F)=\overline{\bigcup\limits_{i=1}^kA_i}.$$

\vskip 2mm

\no{\bf Case $3$} \ There is an integer $s$ such that $f_i^T=1 \
f_i^I=f_i^F=0$, $1\leq i\leq s$ but $f_j^T=f_j^I=0, f_j^F=1$ for
$s+1\leq j\leq k$.

\vskip 2mm

In this case,

$$\bigcup\limits_{i=1}^k(A_i,f_i)=\bigcup\limits_{i=1}^sA_i\bigcup
\overline{\bigcup\limits_{i=s+1}^kA_i}.$$

\vskip 2mm

\no{\bf Case $4$} \ There is an integer $l$ such that
$f_l^T\not=1$ or $f_l^F\not=1$.

\vskip 2mm

In this case, the union is a general neutrosophic set. It can not
be represented by abstract sets.

If $A\bigcap B=\emptyset$, define the function value of a function
$f$ on the union set $A\bigcup B$ to be

$$f(A\bigcup B)=f(A)+f(B)$$

\no and

$$f(A\bigcap B)=f(A)f(B)$$.

\no Then if $A\bigcap B\not=\emptyset$, we get that

$$f(A\bigcup B)= f(A)+f(B)- f(A)f(B).$$

 Generally, by applying the Inclusion-Exclusion Principle to a union of sets,
we get the following formulae.

$$f(\bigcap\limits_{i=1}^lA_i)=\prod\limits_{i=1}^lf(A_i),$$

$$f(\bigcup\limits_{i=1}^kA_i)=\sum\limits_{j=1}^k(-1)^{j-1}\prod\limits_{s=1}^jf(A_s).$$

\vskip 6mm

\no{\bf \S $1.2$ \ Algebraic Structures}

\vskip 4mm

\no In this section, we recall some conceptions and results
without proofs in algebra, such as, these groups, rings, fields,
vectors $\cdots$, all of these can be viewed as a sole-space
system.

\vskip 4mm

\no{\bf $1.2.1.$ Groups}

\vskip 3mm

\no A set $G$ with a binary operation ¡°$\circ$¡±, denoted by
$(G;\circ )$, is called a {\it group} if $x\circ y\in G$ for
$\forall x,y\in G$ such that the following conditions hold.

\vskip 3mm

($i$) $(x\circ y)\circ z=x\circ (y\circ z)$ for $\forall x,y,z\in
G$;

($ii$) There is an element $1_G, 1_G\in G$ such that $x\circ
1_G=x$;

($iii$) For $\forall x\in G$, there is an element $y, y\in G$,
such that $x\circ y=1_G$.

\vskip 2mm

A group $G$ is {\it abelian} if the following additional condition
holds. \vskip 2mm

($iv$) For $\forall x,y\in G$, $x\circ y=y\circ x.$

A set $G$ with a binary operation ¡°$\circ$¡± satisfying the
condition $(i)$ is called a {\it semigroup}. Similarly, if it
satisfies the conditions ($i$) and ($iv$), then it is called a
{\it abelian semigroup}.

Some examples of groups are as follows.

\vskip 3mm

($1$) $(R \ ;+)$ and $(R \ ;\cdot )$, where $R$ is the set of real
numbers.

($2$) $(U_2;\cdot )$, where $U_2=\{1,-1\}$ and generally,
$(U_n;\cdot )$, where $U_n=\{e^{i\frac{2\pi k}{n}}, k=1,2,\cdots
,n\}$.

($3$) For a finite set $X$, the set $Sym X$ of all permutations on
$X$ with respect to permutation composition.

The cases ($1$) and ($2$) are abelian group, but ($3$) is not in
general.

A subset $H$ of a group $G$ is said to be {\it subgroup} if $H$ is
also a group under the same operation in $G$, denoted by $H\prec
G$. The following results are well-known.

\vskip 4mm

\no{\bf Theorem $1.2.1$} \ {\it A non-empty subset $H$ of a group
$(G \ ;\circ )$ is a group if and only if for $\forall x,y\in H$,
$x\circ y\in H$.}

\vskip 3mm

\no{\bf Theorem $1.2.2$}(Lagrange theorem) \ {\it For any subgroup
$H$ of a finite group $G$, the order $|H|$ is a divisor of $|G|$.}

\vskip 3mm

For $\forall x\in G$, denote the set $\{xh|\forall h\in H\}$ by
$xH$ and $\{hx|\forall h\in H\}$ by $Hx$. A subgroup $H$ of a
group $(G \ ;\circ )$ is {\it normal}, denoted by $H\triangleleft
G$, if for $\forall x\in G$, $xH=Hx$.

For two subsets $A,B$ of a group $(G \ ;\circ )$, define their
product $A\circ B$ by

$$A\circ B=\{ a\circ b | \ \forall a\in A, \ \forall b\in b \ \}.$$

For a subgroup $H, H\triangleleft G$, it can be shown that

$$(xH)\circ (yH)= (x\circ y)H \ {\rm and} \ (Hx)\circ (Hy)= H(x\circ y).$$

\no  for $\forall x,y\in G$. Whence, the operation "$\circ$" is
closed in the sets $\{xH| x\in G\}=\{Hx|x\in G\}$, denote this set
by $G/H$. We know $G/H$ is also a group by the facts

$$(xH\circ yH)\circ zH=xH\circ (yH\circ zH), \  \forall x,y,z\in G$$

\no{and}

$$(xH)\circ H=xH, \ (xH)\circ (x^{-1}H)=H.$$

For two groups $G,G'$, let $\sigma$ be a mapping from $G$ to $G'$.
If

$$\sigma (x\circ y)=\sigma (x)\circ \sigma (y),$$

\no for $\forall x,y\in G$, then call $\sigma$ a {\it
homomorphism} from $G$ to $G'$. The {\it image} $Im\sigma$ and the
{\it kernel} $Ker\sigma$ of a homomorphism $\sigma: G\rightarrow
G'$ are defined as follows:

$$Im \sigma= G^{\sigma}=\{ \sigma (x) | \ \forall x\in G \ \},$$

$$Ker \sigma=\{ x | \ \forall x\in G, \ \sigma (x)=1_{G'} \ \}.$$

A one to one homomorphism is called a {\it monomorphism} and an
onto homomorphism an {\it epimorphism}. A homomorphism is called a
{\it bijection} if it is one to one and onto. Two groups $G,G'$
are said to be {\it isomorphic} if there exists a bijective
homomorphism $\sigma$ between them, denoted by $G\cong G'$.

\vskip 4mm

\no{\bf Theorem $1.2.3$} \ {\it Let $\sigma: G\rightarrow G'$ be a
homomorphism of group. Then}

$$(G,\circ )/Ker \sigma\cong Im\sigma.$$

\vskip 6mm

\no{\bf $1.2.2.$ Rings}

\vskip 4mm

\no A set $R$ with two binary operations ¡°$+$¡± and ¡°$\circ$¡±,
denoted by $(R \ ;+,\circ )$, is said to be a {\it ring} if
$x+y\in R$, $x\circ y\in R$ for $\forall x,y\in R$ such that the
following conditions hold.

\vskip 3mm

($i$) \ $(R \ ;+)$ is an abelian group;

($ii$) \ $(R \ ;\circ )$ is a semigroup;

($iii$) \ For $\forall x,y,z\in R$, $x\circ (y+z)=x\circ y+x\circ
z$ and $(x+y)\circ z= x\circ z+y\circ z.$

\vskip 2mm

Some examples of rings are as follows.

\vskip 3mm

($1$) \ $(Z \ ;+,\cdot )$, where $Z$ is the set of integers.

($2$) \ $(pZ \ ;+,\cdot )$, where $p$ is a prime number and
$pZ=\{pn| n\in Z\}$.

($3$) \ $({\mathcal M}_n(Z) \ ; +,\cdot )$, where  ${\mathcal
M}_n(Z)$ is a set of $n\times n$ matrices with each entry being an
integer, $n\geq 2$.\vskip 2mm

For a ring $(R \ ;+,\circ )$, if $x\circ y= y\circ x$ for $\forall
x,y\in R$, then it is called a {\it commutative ring}. The
examples of ($1$) and ($2$) are commutative, but ($3$) is not.

If $R$ contains an element $1_R$ such that for $\forall x\in R$,
$x\circ 1_R=1_R\circ x=x$, we call $R$ a {\it ring with unit}. All
of these examples of rings in the above are rings with unit.  For
($1$), the unit is $1$, ($2$) is $Z$ and ($3$) is $I_{n\times n}$.

The unit of $(R \ ;+)$ in a ring $(R \ ;+,\circ )$ is called {\it
zero}, denoted by $0$. For $\forall a, b\in R$, if

$$a\circ b \ = \ 0,$$

\no then $a$ and $b$ are called {\it divisors of zero}. In some
rings, such as the $(Z \ ;+,\cdot )$ and $(pZ \ ;+,\cdot )$, there
must be $a$ or $b$ be $0$. We call it only has a {\it trivial
divisor of zero}. But in the ring $(pqZ \ ;+,\cdot )$ with $p,q$
both being prime, since

$$pZ \  \cdot \ qZ \ = \ 0$$

\no and $pZ\not=0$, $qZ\not= 0$, we get non-zero divisors of zero,
which is called to have {\it non-trivial divisors of zero}. The
ring $({\mathcal M}_n(Z); +,\cdot )$ also has non-trivial divisors
of zero, since

\[\left(
\begin{array}{llcl}
1 & 1 & \cdots & 1\\
0 & 0 & \cdots & 0\\
\vdots & \vdots & \ddots & \vdots \\
0 & 0 & \cdots & 0\\
\end{array}
\right) \cdot \left(
\begin{array}{llcl}
0 & 0 & \cdots & 0\\
0 & 0 & \cdots & 0\\
\vdots & \vdots & \ddots & \vdots \\
1 & 1 & \cdots & 1\\
\end{array}
\right) = O_{n\times n}. \]

A {\it division ring} is a ring which has no non-trivial divisors
of zero and an {\it integral domain} is a commutative ring having
no non-trivial divisors of zero.

A {\it body} is a ring $(R \ ;+,\circ )$ with a unit, $|R|\geq 2$
and $(R\setminus \{0\};\circ )$ is a group and a {\it field} is a
commutative body. The examples ($1$) and ($2$) of rings are
fields. The following result is well-known.

\vskip 4mm

\no{\bf Theorem $1.2.4$} \ {\it Any finite integral domain is a
field.}

\vskip 3mm

A non-empty subset $R'$ of a ring $(R \ ;+,\circ )$ is called a
{\it subring} if $(R' \ ;+,\circ )$ is also a ring. The following
result for subrings can be obtained immediately by definition.

\vskip 4mm

\no{\bf Theorem $1.2.5$} \ {\it For a subset $R'$ of a ring $(R \
;+,\circ )$, if

($i$) $(R' \ ;+)$ is a subgroup of $(R \ ;+)$,

($ii$) $R'$ is closed under the operation ¡°$\circ$¡±,

\no then $(R' \ ;+,\circ )$ is a subring of $(R \ ,+.\circ )$.}

\vskip 3mm

An {\it ideal} $I$ of a ring $(R \ ;+,\circ )$ is a non-void
subset of R with properties:

($i$) \ $(I \ ;+)$ is a subgroup of $(R \ ;+)$;

($ii$) \ $a\circ x\in I$ and $x\circ a\in I$ for $\forall a\in I,
\forall x\in R$.

Let $(R \ ;+,\circ)$ be a ring. A chain

$$R\succ R_1\succ\cdots\succ  R_l=\{1_{\circ}\}$$

\no satisfying that $R_{i+1}$ is an ideal of $R_i$ for any integer
$i, 1\leq i\leq l$, is called an {\it ideal chain} of $(R \
,+,\circ)$. A ring whose every ideal chain only has finite terms
is called an {\it Artin ring}. Similar to normal subgroups,
consider the set $x+I$ in the group $(R \ ;+)$. Calculation shows
that $R/I= \{x+I| \ x\in R\}$ is also a ring under these
operations ¡°$+$¡± and ¡°$\circ$¡±. Call it a {\it quotient ring
of $R$ to $I$}.

For two rings $(R \ ;+,\circ ),(R' \ ;\ast ,\bullet )$, let
$\iota$ be a mapping from $R$ to $R'$. If

$$\iota (x+y)=\iota (x)\ast\iota (y),$$

$$\iota (x\circ y)=\iota (x)\bullet \iota (y),$$

\no for $\forall x,y\in R$, then $\iota$ is called a {\it
homomorphism} from $(R \ ;+,\circ )$ to $(R' \ ;\ast ,\bullet )$.
Similar to Theorem $2.3$, we know that

\vskip 4mm

\no{\bf Theorem $1.2.6$} \ {\it Let $\iota: R\rightarrow R'$ be a
homomorphism from $(R \ ;+,\circ )$ to $(R' \ ;\ast ,\bullet )$.
Then}

$$(R \ ;+,\circ )/Ker \iota\cong Im\iota.$$

\vskip 6mm

\no{\bf $1.2.3$ Vector spaces}

\vskip 4mm

\no A {\it vector space} or {\it linear space} consists of the
following:

($i$) a field $F$ of scalars;

($ii$) a set $V$ of objects, called vectors;

($iii$) an operation, called vector addition, which associates
with each pair of vectors ${\bf a,b}$ in $V$ a vector ${\bf a+ b}$
in $V$, called the sum of ${\bf a}$ and ${\bf b}$, in such a way
that

($1$) addition is commutative, ${\bf a+ b = b+ a}$;

($2$) addition is associative, ${\bf (a+ b)+ c= a+( b+ c)}$;

($3$) there is a unique vector ${\bf 0}$ in $V$, called the zero
vector, such that ${\bf a+ 0= a}$ for all ${\bf a}$ in $V$;

($4$) for each vector ${\bf a}$ in $V$there is a unique vector
${\bf -a}$ in $V$ such that ${\bf a+(-a)= 0}$;

($iv$) an operation ¡°$\cdot$¡±, called scalar multiplication,
which associates with each scalar $k$ in $F$ and a vector ${\bf
a}$ in $V$ a vector $k\cdot {\bf a}$ in $V$, called the product of
$k$ with ${\bf a}$, in such a way that

($1$) $1\cdot {\bf a}={\bf a}$ for every ${\bf a}$ in $V$;

($2$) $(k_1k_2)\cdot {\bf a}=k_1(k_2\cdot {\bf a});$

($3$) $k\cdot ({\bf a}+{\bf b})=k\cdot {\bf a}+k\cdot {\bf b};$

($4$) $(k_1+k_2)\cdot {\bf a}=k_1\cdot {\bf a}+k_2\cdot {\bf a}.$

\no We say that $V$ is a {\it vector space over the field $F$,
denoted by $(V \ ;+,\cdot )$.}

Some examples of vector spaces are as follows.

($1$) {\it The $n$-tuple space $R^n$ over the real number field
$R$}. \ Let $V$ be the set of all $n$-tuples $(x_1,x_2,\cdots
,x_n)$ with $x_i\in R, 1\leq i\leq n$. If $\forall {\bf
a}=(x_1,x_2,\cdots ,x_n)$, ${\bf b}=(y_1,y_2,\cdots , y_n)\in V$,
then the sum of ${\bf a}$ and ${\bf b}$ is defined by

$${\bf a}+{\bf b}=(x_1+y_1,x_2+y_2,\cdots ,x_n+y_n).$$

\no The product of a real number $k$ with ${\bf a}$ is defined by

$$k{\bf a}=(kx_1,kx_2,\cdots , kx_n).$$

($2$) {\it The space $Q^{m\times n}$ of $m\times n$ matrices over
the rational number field $Q$.} \ Let $Q^{m\times n}$ be the set
of all $m\times n$ matrices over the natural number field $Q$. The
sum of two vectors $A$ and $B$ in $Q^{m\times n}$ is defined by

$$(A+B)_{ij}=A_{ij}+B_{ij},$$

\no and the product of a rational number $p$ with a matrix $A$ is
defined by

$$(pA)_{ij}=pA_{ij}.$$

A {\it subspace} $W$ of a vector space $V$ is a subset $W$ of $V$
which is itself a vector space over $F$ with the operations of
vector addition and scalar multiplication on $V$. The following
result for subspaces is known in references $[6]$ and $[33]$.

\vskip 4mm

\no{\bf Theorem $1.2.7$} \ {\it A non-empty subset $W$ of a vector
space $(V \ ;+,\cdot )$ over the field $F$ is a subspace of $(V \
;+,\cdot )$ if and only if for each pair of vectors ${\bf a},{\bf
b}$ in $W$ and each scalar $k$ in $F$ the vector $k\cdot {\bf
a}+{\bf b}$ is also in $W$.}

\vskip 3mm

Therefore, the intersection of two subspaces of a vector space $V$
is still a subspace of $V$. Let $U$ be a set of some vectors in a
vector space $V$ over $F$. The subspace spanned by $U$ is defined
by

$$\left< U \right> \ = \{ \ k_1\cdot {\bf a_1}+k_2\cdot
{\bf a_2}+\cdots +k_l\cdot {\bf a}_l \ | \ l\geq 1, k_i\in F,
  \ {\rm and} \ {\bf a_j}\in S, 1\leq i\leq
l \ \}.$$

A subset $W$ of $V$ is said to be linearly dependent if there
exist distinct vectors ${\bf a_1,a_2,\cdots ,a_n}$ in $W$ and
scalars $k_1,k_2,\cdots ,k_n$ in $F$, not all of which are $0$,
such that

$$k_1\cdot {\bf a_1}+k_2\cdot {\bf a_2}+\cdots +k_n\cdot {\bf a_n}= {\bf 0}.$$

For a vector space $V$, its {\it basis} is a linearly independent
set of vectors in $V$ which spans the space $V$. Call a space $V$
{\it finite-dimensional} if it has a finite basis. Denoted by $dim
V$ the number of elements in a basis of $V$.

For two subspaces $U,W$ of a space $V$, the sum of subspaces $U,
W$ is defined by

$$U+W=\{ \ {\bf u}+{\bf w} \ | \ {\bf u}\in U, \ {\bf w}\in W \ \}.$$

Then, we have results in the following ([6][33]).

\vskip 4mm

\no{\bf Theorem $1.2.8$} \ {\it Any finite-dimensional vector
space $V$ over a field $F$ is isomorphic to one and only one space
$F^n$, where $n=dimV$.}

\vskip 3mm

\no{\bf Theorem $1.2.9$} \ {\it If $W_1$ and $W_2$ are
finite-dimensional subspaces of a vector space $V$, then $W_1+W_2$
is finite-dimensional and}

$$dimW_1+dimW_2 = dim(W_1\bigcap W_2)+dim(W_1+W_2).$$

\vskip 6mm

\no{\bf \S $1.3$ \ Algebraic Multi-Spaces}

\vskip 4mm

\no The notion of a multi-space was introduced by Smarandache in
1969 ([$86$]). Algebraic multi-spaces had be researched in
references $[58]-[61]$ and $[103]$. Vasantha Kandasamy researched
various bispaces in $[101]$, such as those of bigroups,
bisemigroups, biquasigroups, biloops, bigroupoids, birings,
bisemirings, bivectors, bisemivectors, bilnear-rings, $\cdots$,
etc., considered two operation systems on two different sets.

\vskip 4mm

\no{\bf $1.3.1.$ Algebraic multi-spaces}

\vskip 3mm

\no{\bf Definition $1.3.1$} \ {\it For any integers $n,i$, $n\geq
2$ and $1\leq i\leq n$, let $A_i$ be a set with ensemble of law
$L_i$, and the intersection of $k$ sets $A_{i_1},A_{i_2},\cdots ,
A_{i_k}$ of them constrains the law $I(A_{i_1},A_{i_2},\cdots ,
A_{i_k})$. Then the union $\widetilde{A}$

$$\widetilde{A} \ = \ \bigcup\limits_{i=1}^n A_i$$

\no is called a multi-space.}

Notice that in this definition, each law may be contain more than
one binary operation. For a binary operation ¡°$\times$¡±, if
there exists an element $1_{\times}^l$ (or $1_{\times}^r$) such
that

$$1_{\times}^l\times a = a \ \ {\rm or} \ \ a\times 1_{\times}^r=a$$

\no for $\forall a\in A_i, 1\leq i\leq n$, then $1_{\times}^l$
($1_{\times}^r$) is called a {\it left (right) unit}. If
$1_{\times}^l$ and $1_{\times}^r$ exist simultaneously, then there
must be

$$1_{\times}^l=1_{\times}^l\times 1_{\times}^r=1_{\times}^r=1_{\times}.$$

\no Call $1_{\times}$ a {\it unit} of $A_i$.

\vskip 3mm

\no{\bf Remark $1.3.1$} \ In Definition $1.3.1$, the following
three cases are permitted:

($i$) \ $A_1=A_2=\cdots = A_n$, i.e., $n$ laws on one set.

($ii$) \ $L_1=L_2=\cdots =L_n$, i.e., n set with one law

($iii$) \ there exist integers $s_1,s_2,\cdots ,s_l$ such that
$I(s_j)=\emptyset , 1\leq j\leq l$, i.e., some laws on the
intersections may be not existed.

\vskip 2mm

We give some examples for Definition $1.3.1$.

\vskip 3mm

\no{\bf Example $1.3.1$} \ Take $n$ disjoint two by two cyclic
groups $C_1,C_2,\cdots ,C_n, \ n\geq 2$ with

$$C_1= (\left< a\right>;+_1), C_2= (\left< b\right>;+_2),\cdots ,
C_n= (\left< c\right>;+_n).$$

\no Where ¡°$+_1,+_2,\cdots , +_n$¡± are $n$ binary operations.
Then their union

$$\widetilde{C}= \bigcup\limits_{i=1}^n C_i$$

\no is a multi-space with the empty intersection laws. In this
multi-space, for $\forall x,y\in \widetilde{C}$, if $x,y\in C_k$
for some integer $k$, then we know $x+_ky\in C_k$. But if $x\in
C_s$, $y\in C_t$ and $s\not=t$, then we do not know which binary
operation between them and what is the resulting element
corresponds to them.

A general multi-space of this kind is constructed by choosing $n$
algebraic systems $A_1,A_2, \cdots , A_n$ satisfying that

$$A_i\bigcap A_j=\emptyset \ \ {\rm and} \ \ O(A_i)\bigcap O(A_j)=\emptyset,$$

\no for any integers $i,j, i\not=j, \ 1\leq i, j\leq n$, where
$O(A_i)$ denotes the binary operation set in $A_i$. Then

$$\widetilde{A}=\bigcup\limits_{i=1}^nA_i$$

\no with $O(\widetilde{A})=\bigcup\limits_{i=1}^nO(A_i)$ is a
multi-space. This kind of multi-spaces can be seen as a model of
spaces with a empty intersection.

\vskip 3mm

\no{\bf Example $1.3.2$} \ Let $(G \ ;\circ )$ be a group with a
binary operation ¡°$\circ$¡±. Choose $n$ different elements
$h_1,h_2,\cdots ,h_n$, $n\geq 2$ and make the extension of the
group $(G \ ;\circ )$ by $h_1,h_2,\cdots ,h_n$ respectively as
follows:

$(G\bigcup\{h_1\}; \times_1)$, \ where the binary operation
$\times_1=\circ$ for elements in $G$, otherwise, new operation;

$(G\bigcup\{h_2\}; \times_2)$, \ where the binary operation
$\times_2=\circ$ for elements in $G$, otherwise, new operation;

$\cdots\cdots\cdots\cdots\cdots\cdots ;$

$(G\bigcup\{h_n\}; \times_n)$, \ where the binary operation
$\times_n=\circ$ for elements in $G$, otherwise, new operation.

Define

$$\widetilde{G}= \bigcup\limits_{i=1}^n(G\bigcup\{h_i\}; \times_i).$$

\no Then $\widetilde{G}$ is a multi-space with binary operations
¡°$\times_1,\times_2, \cdots , \times_n$¡±. In this multi-space,
for $\forall x,y\in\widetilde{G}$, unless the exception cases
$x=h_i,y=h_j$ and $i\not=j$, we know the binary operation between
$x$ and $y$ and the resulting element by them.

For $n=3$, this multi-space can be shown as in Fig.$1.2$, in where
the central circle represents the group $G$ and each angle field
the extension of $G$. Whence, we call this kind of multi-space a
{\it fan multi-space}. \vskip 2mm

\includegraphics[bb=-90 10 200 110]{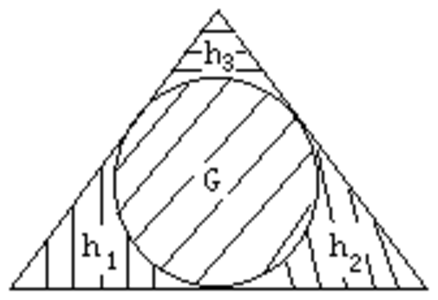}

\vskip 2mm

\c{\bf Fig.$1.2$}

\vskip 2mm

Similarly, we can also use a ring $R$ to get fan multi-spaces. For
example, let $(R \ ;+,\circ )$ be a ring and let $r_1,r_2,\cdots
,r_s$ be two by two different elements. Make these extensions of
$(R \ ;+,\circ )$ by $r_1,r_2,\cdots ,r_s$ respectively as
follows:

$(R\bigcup\{r_1\}; +_1,\times_1)$, \ where  binary operations
$+_1=+$, \ $\times_1=\circ$ for elements in $R$, otherwise, new
operation;

$(R\bigcup\{r_2\}; +_2, \times_2)$, \ where  binary operations
$+_2=+$, \ $\times_2=\circ$ for elements in $R$, otherwise, new
operation;

$\cdots\cdots\cdots\cdots\cdots\cdots ;$

$(R\bigcup\{r_s\}; +_s, \times_s)$, \ where  binary operations
$+_s=+$, \ $\times_s=\circ$ for elements in $R$, otherwise, new
operation.

Define

$$\widetilde{R}=\bigcup\limits_{j=1}^s(R\bigcup\{r_j\}; +_j,
\times_j).$$

\no Then $\widetilde{R}$ is a fan multi-space with ring-like
structure. Also we can define a fan multi-space with field-like,
vector-like, semigroup-like,$\cdots$, etc. structures.

These multi-spaces constructed in Examples $1.3.1$ and $1.3.2$ are
not {\it completed}, i.e., there exist some elements in this space
not have binary operation between them. In algebra, we wish to
construct a {\it completed multi-space}, i.e., there is a binary
operation between any two elements at least and their resulting is
still in this space. The following example is a completed
multi-space constructed by applying {\it Latin squares} in the
combinatorial design.

\vskip 3mm

\no{\bf Example $1.3.3$} \ Let $S$ be a finite set with $|S|=n\geq
2$. Constructing an $n\times n$ Latin square by elements in $S$,
i.e., every element just appears one time on its each row and each
column. Now choose $k$ Latin squares $M_1,M_2,\cdots , M_k$,
$k\leq\prod\limits_{s=1}^ns!$.

By a result in the reference $[83]$, there are at least
$\prod\limits_{s=1}^ns!$ distinct $n\times n$ Latin squares.
Whence, we can always choose $M_1,M_2,\cdots , M_k$ distinct two
by two. For a Latin square $M_i, 1\leq i\leq k$, define an
operation ¡°$\times_i$¡±as follows:

$$\times_i: (s,f)\in S\times S \rightarrow (M_i)_{sf}.$$

The case of $n=3$ is explained in the following. Here
$S=\{1,2,3\}$ and there are $2$ Latin squares $L_1,L_2$ as
follows:

\[
L_1=\left(
\begin{array}{lll}
1 & 2 & 3\\
2 & 3 & 1\\
3 & 1 & 2\\
\end{array}
\right) \quad\quad L_2= \left(
\begin{array}{lll}
1 & 2 & 3\\
3 & 1 & 2\\
2 & 3 & 1\\
\end{array}
\right)  \].

Therefore, by the Latin square $L_1$, we get an operation
¡°$\times_1$¡±as in table $1.3.1$.\vskip 2mm

\begin{center}
\begin{tabular}{c|ccc}
$\times_1$ & \ $1$ \ & \ $2$ \ & \ $3$ \ \\ \hline
$1$ & \ $1$ \ & \ $2$ \ & \ $3$ \ \\
$2$ & \ $2$ \ & \ $3$ \ & \ $1$ \ \\
$3$ & \ $3$ \ & \ $1$ \ & \ $2$ \ \\
\end{tabular}
\end{center}
\vskip 2mm

\c{\bf table $1.3.1$}

\vskip 2mm

\no and by the Latin square $L_2$, we also get an operation
¡°$\times_2$¡± as in table $1.3.2$.\vskip 2mm

\begin{center}
\begin{tabular}{c|ccc}
$\times_2$ & \ $1$ \ & \ $2$ \ & \ $3$ \ \\ \hline
$1$ & \ $1$ \ & \ $2$ \ & \ $3$ \ \\
$2$ & \ $3$ \ & \ $1$ \ & \ $2$ \ \\
$3$ & \ $2$ \ & \ $3$ \ & \ $1$ \ \\
\end{tabular}
\end{center}

\vskip 2mm

\c{\bf table $1.3.2$}

\vskip 2mm

For $\forall x,y,z\in S$ and two operations ¡°$\times_i$¡± and
¡°$\times_j$¡±, $1\leq i,j\leq k$, define

$$x\times_iy\times_jz = (x\times_iy)\times_jz.$$

\no For example, in the case $n=3$, we know that

$$1\times_12\times_23 = (1\times_2)\times_23 = 2\times_23 = 2;$$

\no and

$$2\times_13\times_22 = (2\times_13)\times_22=1\times_23 = 3.$$

\no Whence $S$ is a completed multi-space with $k$ operations.

The following example is also a completed multi-space constructed
by an algebraic system.

\vskip 3mm

\no{\bf Example $1.3.4$} \ For constructing a completed
multi-space, let $(S \ ;\circ )$ be an algebraic system, i.e.,
$a\circ b\in S$ for $\forall a,b\in S$. Whence, we can take $C,
C\subseteq S$ being a cyclic group. Now consider a partition of
$S$

$$S \ = \ \bigcup\limits_{k=1}^mG_k$$

\no with $m\geq 2$ such that $G_i\bigcap G_j=C$ for $\forall i,j,
1\leq i,j \leq m$.

For an integer $k, 1\leq k\leq m$, assume
$G_k=\{g_{k1},g_{k2},\cdots , g_{kl}\}$. We define an operation
¡°$\times_k$¡±on $G_k$ as follows, which enables $(G_k;\times_k)$
to be a cyclic group.

$$g_{k1}\times_kg_{k1}=g_{k2},$$

$$g_{k2}\times_kg_{k1}=g_{k3},$$

$$\cdots\cdots\cdots\cdots\cdots\cdots ,$$

$$g_{k(l-1)}\times_kg_{k1}=g_{kl},$$

\no and

$$g_{kl)}\times_kg_{k1}=g_{k1}.$$

\no Then $S \ = \ \bigcup\limits_{k=1}^mG_k$ is a completed
multi-space with $m+1$ operations.

The approach used in Example $1.3.4$ enables us to construct a
complete multi-spaces $\widetilde{A}=\bigcup\limits_{i=1}^n$ with
$k$  operations for $k\geq n+1$, i.e., the intersection law
$I(A_1,A_2,\cdots , A_n)\not=\emptyset$.

\vskip 4mm

\no{\bf Definition $1.3.2$} \ {\it A mapping $f$ on a set $X$ is
called faithful if $f(x)=x$ for $\forall x\in X$, then $f=1_X$,
the unit mapping on $X$ fixing each element in $X$. }

\vskip 3mm

Notice that if $f$ is faithful and $f_1(x)=f(x)$ for $\forall x\in
X$, then $f_1^{-1}f=1_X$, i.e., $f_1=f$.

For each operation ¡°$\times$¡± and a chosen element $g$ in a
subspace $A_i, A_i\subset\widetilde{A}, 1\leq i\leq n$, there is a
{\it left-mapping} $f_g^l: A_i\rightarrow A_i$ defined by

$$f_g^l: a\rightarrow g\times a, \ \ a\in A_i.$$

\no Similarly, we can also define the {\it right-mapping} $f_g^r$.

We adopt the following convention for multi-spaces in this book.

\vskip 4mm

\no{\bf Convention $1.3.1$} \ {\it Each operation ¡°$\times$¡±in a
subset $A_i, A_i\subset\widetilde{A}, 1\leq i\leq n$ is faithful,
i.e., for $\forall g\in A_i$, $\varsigma : g\rightarrow f_g^l$ ( \
or $\tau : g\rightarrow f_g^r$ ) is faithful.}

\vskip 3mm

Define the kernel $Ker\varsigma$ of a mapping $\varsigma$ by

$${\rm Ker}\varsigma =\{g| g\in A_i \ {\rm and} \ \varsigma (g)= 1_{A_i}\}.$$

\no Then Convention $1.3.1$ is equivalent to the next convention.

\vskip 4mm

\no{\bf Convention $1.3.2$} \ {\it For each $\varsigma :
g\rightarrow f_g^l$ ( \ or $\varsigma : g\rightarrow f_g^r$ )
induced by an operation ¡°$\times$¡± has kernel}

$${\rm Ker}\varsigma =\{1_{\times}^l\}$$

\no{\it if $1_{\times}^l$ exists. Otherwise, ${\rm Ker}\varsigma
=\emptyset$.}

\vskip 3mm

We have the following results for multi-spaces $\widetilde{A}$.

\vskip 4mm

\no{\bf Theorem $1.3.1$} \ {\it For a multi-space $\widetilde{A}$
and an operation ¡°$\times$¡±, the left unit $1_{\times}^l$ and
right unit $1_{\times}^r$ are unique if they exist.}

\vskip 3mm

{\it Proof} \ If there are two left units $1_{\times}^l,
I_{\times}^l$ in a subset $A_i$ of a multi-space $\widetilde{A}$,
then for $\forall x\in A_i$, their induced left-mappings
$f_{1_{\times}^l}^l$ and $f_{I_{\times}^l}^l$ satisfy

$$f_{1_{\times}^l}^l(x)=1_{\times}^l\times x=x$$

\no and

$$f_{I_{\times}^l}^l(x)=I_{\times}^l\times x=x.$$

Therefore, we get that $f_{1_{\times}^l}^l=f_{I_{\times}^l}^l$.
Since the mappings $\varsigma_1 : 1_{\times}^l\rightarrow
f_{1_{\times}^l}^l$ and $\varsigma_2 : I_{\times}^l\rightarrow
f_{I_{\times}^l}^l$ are faithful, we know that

$$1_{\times}^l = I_{\times}^l.$$

Similarly, we can also prove that the right unit $1_{\times}^r$ is
also unique.\quad\quad $\natural$

For two elements $a, b$ of a multi-space $\widetilde{A}$, if
$a\times b= 1_{\times}^l$, then $b$ is called a {\it left-inverse}
of $a$. If $a\times b= 1_{\times}^r$, then $a$ is called a {\it
right-inverse} of $b$. Certainly, if $a\times b= 1_{\times}$, then
$a$ is called an {\it inverse} of $b$ and $b$ an {\it inverse} of
$a$.

\vskip 4mm

\no{\bf Theorem $1.3.2$} \ {\it For a multi-space $\widetilde{A}$,
$a\in\widetilde{A}$, the left-inverse and right-inverse of $a$ are
unique if they exist.}

\vskip 3mm

{\it Proof} \ Notice that $\kappa_a : x\rightarrow ax$ is
faithful, i.e., ${\rm Ker}\kappa =\{1_{\times}^l\}$ for
$1_{\times}^l$ existing now.

If there exist two left-inverses $b_1, b_2$ in $\widetilde{A}$
such that $a\times b_1= 1_{\times}^l$ and $a\times b_2=
1_{\times}^l$, then we know that

$$b_1 = b_2 = 1_{\times}^l.$$

Similarly, we can also prove that the right-inverse of $a$ is also
unique. \quad\quad $\natural$

\vskip 4mm

\no{\bf Corollary $1.3.1$} \ {\it If ¡°$\times$¡± is an operation
of a multi-space $\widetilde{A}$ with unit $1_{\times}$, then the
equation}

$$a\times x = b $$

\no{\it has at most one solution for the indeterminate $x$.}

\vskip 3mm

{\it Proof} \ According to Theorem $1.3.2$, we know there is at
most one left-inverse $a_1$ of $a$ such that $a_1\times a =
1_{\times}$. Whence, we know that

$$x = a_1\times a\times x = a_1\times b.\quad\quad \natural$$

We also get a consequence for solutions of an equation in a
multi-space by this result.

\vskip 4mm

\no{\bf Corollary $1.3.2$} \ {\it Let $\widetilde{A}$ be a
multi-space with a operation set $O(\widetilde{A})$. Then the
equation}

$$a\circ x = b $$

\no{\it has at most $o(\widetilde{A})$ solutions, where
¡°$\circ$¡±is any binary operation of $\widetilde{A}$.}

\vskip 3mm

Two multi-spaces $\widetilde{A_1}, \widetilde{A_2}$ are said to be
{\it isomorphic} if there is a one to one mapping $\zeta
:\widetilde{A_1}\rightarrow\widetilde{A_2}$ such that for $\forall
x,y\in\widetilde{A}_1$ with binary operation ¡°$\times$¡±, $\zeta
(x), \zeta (y)$ in $\widetilde{A_2}$ with binary operation
¡°$\circ$¡± satisfying the following condition

$$\zeta (x\times y)=\zeta (x)\circ \zeta (y).$$

\no If $\widetilde{A_1}= \widetilde{A_2}=\widetilde{A}$, then an
isomorphism between $\widetilde{A_1}$ and  $\widetilde{A_2}$ is
called an {\it automorphism} of $\widetilde{A}$. All automorphisms
of $\widetilde{A}$ form a group under the composition operation
between mappings, denoted by ${\rm Aut}\widetilde{A}$.

Notice that ${\rm Aut}Z_n\cong Z_n^*$, where $Z_n^*$ is the group
of reduced residue class mod$n$ under the multiply operation (
$[108]$ ). It is known that $|{\rm Aut}Z_n|=\varphi (n)$, where
$\varphi (n)$ is the Euler function. We know the automorphism
group of the multi-space $\widetilde{C}$ in Example $1.3.1$ is

$${\rm Aut}\widetilde{C}=S_n[Z_n^*].$$

\no Whence, $|{\rm Aut}\widetilde{C}|=\varphi (n)^nn!$. For
Example $1.3.3$, determining its automorphism group is a more
interesting  problem for the combinatorial design ( see also the
final section in this chapter).

\vskip 6mm

\no{\bf $1.3.2$ Multi-Groups}

\vskip 4mm

\no The conception of multi-groups is a generalization of
classical algebraic structures, such as those of groups, fields,
bodies, $\cdots$, etc., which is defined in the following
definition.

\vskip 4mm

\no{\bf Definition $1.3.3$} \ {\it Let
$\widetilde{G}=\bigcup\limits_{i=1}^n G_i$ be a complete
multi-space with an operation set $O(\widetilde{G})=\{\times_i,
1\leq i\leq n\}$. If $(G_i;\times_i)$ is a group for any integer
$i, 1\leq i\leq n$ and for $\forall x,y,z\in \widetilde{G}$ and
$\forall\times,\circ\in O(\widetilde{G})$, $\times\not= \circ$,
there is one operation, for example the operation ¡°$\times$¡±
satisfying the distribution law to the operation
¡°$\circ$¡±provided all of these operating results exist , i.e.,

$$x\times (y\circ z) = (x\times y)\circ (x\times z),$$

$$(y\circ z)\times x = (y\times x)\circ (z\times x),$$

\no then $\widetilde{G}$ is called a multi-group.}

\vskip 3mm

\no{\bf Remark $1.3.2$} \ The following special cases for $n=2$
convince us that multi-groups are a generalization of groups,
fields and bodies, $\cdots$, etc..\vskip 2mm

($i$) \ If $G_1=G_2=\widetilde{G}$, then $\widetilde{G}$ is a
body.

($ii$) If $(G_1;\times_1)$ and $(G_2;\times_2)$ are commutative
groups, then $\widetilde{G}$ is a field.\vskip 2mm

For a multi-group $\widetilde{G}$ and a subset
$\widetilde{G_1}\subset\widetilde{G}$, if $\widetilde{G_1}$ is
also a multi-group under a subset $O(\widetilde{G_1}),
O(\widetilde{G_1})\subset O(\widetilde{G})$, then
$\widetilde{G}_1$ is called a {\it sub-multi-group} of
$\widetilde{G}$, denoted by $\widetilde{G_1} \preceq
\widetilde{G}$. We get a criterion for sub-multi-groups in the
following.

\vskip 4mm

\no{\bf Theorem $1.3.3$} \ {\it For a multi-group
$\widetilde{G}=\bigcup\limits_{i=1}^nG_i$ with an operation set
$O(\widetilde{G})=\{\times_i|1\leq i\leq n\}$, a subset
$\widetilde{G_1}\subset\widetilde{G}$ is a sub-multi-group of
$\widetilde{G}$ if and only if $(\widetilde{G_1}\bigcap
G_k;\times_k)$ is a subgroup of $(G_k;\times_k)$ or
$\widetilde{G_1}\bigcap G_k=\emptyset$ for any integer $k, 1\leq
k\leq n$.}

\vskip 3mm

{\it Proof} \ If $\widetilde{G_1}$ is a multi-group with an
operation set $O(\widetilde{G_1})=\{\times_{i_j}|1\leq j\leq
s\}\subset O(\widetilde{G})$, then

$$\widetilde{G_1}=\bigcup\limits_{i=1}^n(\widetilde{G_1}\bigcap
G_i) =\bigcup\limits_{j=1}^sG'_{i_j}$$

\no where $G'_{i_j}\preceq G_{i_j}$ and $(G_{i_j};\times_{i_j})$
is a group. Whence, if $\widetilde{G_1}\bigcap G_k\not=\emptyset$,
then there exist an integer $l, k=i_l$ such that
$\widetilde{G_1}\bigcap G_k=G'_{i_l}$, i.e.,
$(\widetilde{G_1}\bigcap G_k;\times_k)$ is a subgroup of
$(G_k;\times_k)$.

Now if $(\widetilde{G_1}\bigcap G_k;\times_k)$ is a subgroup of
$(G_k;\times_k)$ or $\widetilde{G_1}\bigcap G_k=\emptyset$ for any
integer $k$, let $N$ denote the index set $k$ with
$\widetilde{G_1}\bigcap G_k\not=\emptyset$, then

$$\widetilde{G_1}=\bigcup\limits_{j\in N}(\widetilde{G_1}\bigcap
G_j)$$

\no and $(\widetilde{G_1}\bigcap G_j,\times_j)$ is a group. Since
$\widetilde{G_1}\subset\widetilde{G}$, $O(\widetilde{G_1})\subset
O(\widetilde{G})$, the associative law and distribute law are true
for the $\widetilde{G_1}$. Therefore, $\widetilde{G_1}$ is a
sub-multi-group of $\widetilde{G}$.\quad\quad $\natural$

For finite sub-multi-groups, we get a criterion as in the
following.

\vskip 4mm

\no{\bf Theorem $1.3.4$} \ {\it Let $\widetilde{G}$ be a finite
multi-group with an operation set
$O(\widetilde{G})=\{\times_i|1\leq i\leq n\}$. A subset
$\widetilde{G_1}$ of $\widetilde{G}$ is a sub-multi-group under an
operation subset $O(\widetilde{G_1})\subset O(\widetilde{G})$ if
and only if $(\widetilde{G_1};\times)$ is complete for each
operation ¡°$\times$¡± in $O(\widetilde{G_1})$.}

\vskip 3mm

{\it Proof} \ Notice that for a multi-group $\widetilde{G}$, its
each sub-multi-group $\widetilde{G_1}$ is complete.

Now if $\widetilde{G_1}$ is a complete set under each operation
¡°$\times_i$¡± in $O(\widetilde{G_1})$, we know that
$(\widetilde{G_1}\bigcap G_i; \times_i)$ is a group or an empty
set. Whence, we get that

$$\widetilde{G_1}=\bigcup\limits_{i=1}^n(\widetilde{G_1}\bigcap
G_i).$$

\no Therefore, $\widetilde{G_1}$ is a sub-multi-group of
$\widetilde{G}$ under the operation set $O(\widetilde{G_1})$.
\quad\quad $\natural$

For a sub-multi-group $\widetilde{H}$ of a multi-group
$\widetilde{G}$, $g\in\widetilde{G}$, define

$$g\widetilde{H}=\{g\times h|h\in\widetilde{H}, \times\in O(\widetilde{H})\}.$$

\no Then for $\forall x,y\in\widetilde{G}$,

$$x\widetilde{H}\bigcap y\widetilde{H}=\emptyset  \ \ {\rm or} \ \ x\widetilde{H}=y\widetilde{H}.$$

\no In fact, if $x\widetilde{H}\bigcap
y\widetilde{H}\not=\emptyset$, let $z\in x\widetilde{H}\bigcap
y\widetilde{H}$, then there exist elements
$h_1,h_2\in\widetilde{H}$ and operations ¡°$\times_i$¡± and
¡°$\times_j$¡± such that

$$z = x\times_i h_1 = y\times_j h_2.$$

Since $\widetilde{H}$ is a sub-multi-group, $(\widetilde{H}\bigcap
G_i;\times_i)$ is a subgroup. Whence, there exists an inverse
element $h_1^{-1}$ in $(\widetilde{H}\bigcap G_i;\times_i)$. We
get that

$$x\times_i h_1\times_i h_1^{-1} = y\times_j h_2\times_ih_1^{-1}.$$

\no i.e.,

$$x=y\times_j h_2\times_ih_1^{-1}.$$

\no Whence,

$$x\widetilde{H}\subseteq y\widetilde{H}.$$

Similarly, we can also get that

$$x\widetilde{H}\supseteq y\widetilde{H}.$$

\no Thereafter, we get that

$$x\widetilde{H} = y\widetilde{H}.$$

Denote the union of two set $A$ and $B$ by $A\bigoplus B$ if
$A\bigcap B=\emptyset$. Then the following result is implied in
the previous proof.

\vskip 4mm

\no{\bf Theorem $1.3.5$} \ {\it For any sub-multi-group
$\widetilde{H}$ of a multi-group $\widetilde{G}$, there is a
representation set $T$, $T\subset\widetilde{G}$, such that }

$$\widetilde{G}=\bigoplus\limits_{x\in T}x\widetilde{H}.$$

\vskip 3mm

For the case of finite groups, since there is only one binary
operation ¡°$\times$¡± and $|x\widetilde{H}|=|y\widetilde{H}|$ for
any $x,y\in\widetilde{G}$, We get a consequence in the following,
which is just the Lagrange theorem for finite groups.

\vskip 4mm

\no{\bf Corollary $1.3.3$}(Lagrange theorem) \ {\it For any finite
group $G$, if $H$ is a subgroup of $G$, then $|H|$ is a divisor of
$|G|$.}

\vskip 3mm

For a multi-group $\widetilde{G}$ and $g\in\widetilde{G}$, denote
all the binary operations associative with $g$ by
$\overrightarrow{O(g)}$ and the elements associative with the
binary operation ¡°$\times$¡± by $\widetilde{G}(\times)$. For a
sub-multi-group $\widetilde{H}$ of $\widetilde{G}$, $\times\in
O(\widetilde{H})$, if

$$g\times h\times g^{-1}\in\widetilde{H},$$

\no for $\forall h\in\widetilde{H}$ and $\forall
g\in\widetilde{G}(\times)$, then we call $\widetilde{H}$ a {\it
normal sub-multi-group } of $\widetilde{G}$, denoted by
$\widetilde{H}\triangleleft\widetilde{G}$. If $\widetilde{H}$ is a
normal sub-multi-group of $\widetilde{G}$, similar to the normal
subgroups of groups, it can be shown that
$g\times\widetilde{H}=\widetilde{H}\times g$, where
$g\in\widetilde{G}(\times)$. Thereby we get a result as in the
following.

\vskip 4mm

\no{\bf Theorem $1.3.6$} \ {\it Let
$\widetilde{G}=\bigcup\limits_{i=1}^nG_i$ be a multi-group with an
operation set $O(\widetilde{G})=\{\times_i|1\leq i\leq n\}$. Then
a sub-multi-group $\widetilde{H}$ of $\widetilde{G}$ is normal if
and only if $(\widetilde{H}\bigcap G_i;\times_i)$ is a normal
subgroup of $(G_i;\times_i)$ or $\widetilde{H}\bigcap
G_i=\emptyset$ for any integer $i, 1\leq i\leq n$.}

\vskip 3mm

{\it Proof} \ We have known that

$$\widetilde{H} = \bigcup\limits_{i=1}^n(\widetilde{H}\bigcap G_i).$$

If $(\widetilde{H}\bigcap G_i;\times_i)$ is a normal subgroup of
$(G_i;\times_i)$ for any integer $i, 1\leq i\leq n$, then we know
that

$$g\times_i(\widetilde{H}\bigcap G_i)\times_ig^{-1} = \widetilde{H}\bigcap G_i$$

\no for $\forall g\in G_i, 1\leq i\leq n$. Whence,

$$g\circ\widetilde{H}\circ g^{-1} = \widetilde{H}$$

\no for $\forall \circ\in O(\widetilde{H})$ and $\forall
g\in\overrightarrow{\widetilde{G}(\circ )}$. That is,
$\widetilde{H}$ is a normal sub-multi-group of $\widetilde{G}$.

Now if $\widetilde{H}$ is a normal sub-multi-group of
$\widetilde{G}$, by definition we know that

$$g\circ\widetilde{H}\circ g^{-1} = \widetilde{H}$$

\no for $\forall \circ\in O(\widetilde{H})$ and $\forall
g\in\widetilde{G}(\circ )$. Not loss of generality, we assume that
$\circ = \times_k$, then we get

$$g\times_k(\widetilde{H}\bigcap G_k)\times_k g^{-1}= \widetilde{H}\bigcap G_k.$$

\no Therefore, $(\widetilde{H}\bigcap G_k; \times_k)$ is a normal
subgroup of $(G_k,\times_k)$. Since the operation ¡°$\circ$¡± is
chosen arbitrarily, we know that $(\widetilde{H}\bigcap
G_i;\times_i)$ is a normal subgroup of $(G_i;\times_i)$ or an
empty set for any integer $i$, $1\leq i\leq n$.\quad\quad
$\natural$

For a multi-group $\widetilde{G}$ with an operation set
$O(\widetilde{G})=\{\times_i| \ 1\leq i\leq n\}$, an order of
operations in $O(\widetilde{G})$ is said to be an {\it oriented
operation sequence}, denoted by
$\overrightarrow{O}(\widetilde{G})$. For example, if
$O(\widetilde{G})=\{\times_1,\times_2\times_3\}$, then
$\times_1\succ\times_2\succ\times_3$ is an oriented operation
sequence and $\times_2\succ\times_1\succ\times_3$ is also an
oriented operation sequence.

For a given oriented operation sequence
$\overrightarrow{O}(\widetilde{G})$, we construct a series of
normal sub-multi-group

$$\widetilde{G}\triangleright\widetilde{G}_1\triangleright\widetilde{G}_2
\triangleright\cdots\triangleright\widetilde{G}_m=\{1_{\times_n}\}$$

\no by the following programming.

\vskip 3mm

\no{\bf STEP $1$}: {\it Construct a series

$$\widetilde{G}\triangleright\widetilde{G}_{11}\triangleright\widetilde{G}_{12}
\triangleright\cdots\triangleright\widetilde{G}_{1l_1}$$

\no under the operation ¡°$\times_1$¡±.}

\vskip 2mm

\no{\bf STEP $2$}: {\it If a series

$$\widetilde{G}_{(k-1)l_1}\triangleright\widetilde{G}_{k1}\triangleright\widetilde{G}_{k2}
\triangleright\cdots\triangleright\widetilde{G}_{kl_k}$$

\no has be constructed under the operation ¡°$\times_k$¡± and
$\widetilde{G}_{kl_k}\not=\{1_{\times_n}\}$, then construct a
series

$$\widetilde{G}_{kl_1}\triangleright\widetilde{G}_{(k+1)1}\triangleright\widetilde{G}_{(k+1)2}
\triangleright\cdots\triangleright\widetilde{G}_{(k+1)l_{k+1}}$$

\no under the operation ¡°$\times_{k+1}$¡±.

This programming is terminated until the series

$$\widetilde{G}_{(n-1)l_1}\triangleright\widetilde{G}_{n1}\triangleright\widetilde{G}_{n2}
\triangleright\cdots\triangleright\widetilde{G}_{nl_n}=\{1_{\times_n}\}$$

\no has be constructed under the operation ¡°$\times_n$¡±.}

The number $m$ is called the {\it length of the series of normal
sub-multi-groups}. Call a series of normal sub-multi-group

$$\widetilde{G}\triangleright\widetilde{G}_1\triangleright\widetilde{G}_2
\triangleright\cdots\triangleright\widetilde{G}_n=\{1_{\times_n}\}$$

\no {\it maximal} if there exists a normal sub-multi-group
$\widetilde{H}$ for any integer $k, s, 1\leq k\leq n, 1\leq s\leq
l_k$ such that

$$\widetilde{G}_{ks}\triangleright\widetilde{H}\triangleright\widetilde{G}_{k(s+1)},$$

\no then $\widetilde{H}=\widetilde{G}_{ks}$ or
$\widetilde{H}=\widetilde{G}_{k(s+1)}$. For a maximal series of
finite normal sub-multi-group, we get a result as in the
following.

\vskip 4mm

\no{\bf Theorem $1.3.7$} \ {\it For a finite multi-group
$\widetilde{G}=\bigcup\limits_{i=1}^n G_i$ and an oriented
operation sequence $\overrightarrow{O}(\widetilde{G})$, the length
of the maximal series of normal sub-multi-group in $\widetilde{G}$
is a constant, only dependent on $\widetilde{G}$ itself.}

\vskip 3mm

{\it Proof} \ The proof is by the induction principle on the
integer $n$.

For $n=1$, the maximal series of normal sub-multi-groups of
$\widetilde{G}$ is just a composition series of a finite group. By
the Jordan-H\"{o}lder theorem (see $[73]$ or $[107]$), we know the
length of a composition series is a constant, only dependent on
$\widetilde{G}$. Whence, the assertion is true in the case of
$n=1$.

Assume that the assertion is true for all cases of $n\leq k$. We
prove it is also true in the case of $n=k+1$. Not loss of
generality, assume the order of those binary operations in
$\overrightarrow{O}(\widetilde{G})$ being
$\times_1\succ\times_2\succ\cdots\succ\times_n$ and the
composition series of the group $(G_1,\times_1)$ being

$$G_1\triangleright G_2\triangleright\cdots\triangleright G_s=\{1_{\times_1}\}.$$

By the Jordan-H\"{o}lder theorem, we know the length of this
composition series is a constant, dependent only on
$(G_1;\times_1)$. According to Theorem $3.6$, we know a maximal
series of normal sub-multi-groups of $\widetilde{G}$ gotten by
STEP $1$ under the operation ¡°$\times_1$¡± is

$$\widetilde{G}\triangleright\widetilde{G}\setminus (G_1\setminus G_2)\triangleright
\widetilde{G}\setminus (G_1\setminus
G_3)\triangleright\cdots\triangleright\widetilde{G}\setminus
(G_1\setminus \{1_{\times_1}\}).$$

Notice that $\widetilde{G}\setminus
(G_1\setminus\{1_{\times_1}\})$ is still a multi-group with less
or equal to $k$ operations. By the induction assumption, we know
the length of the maximal series of normal sub-multi-groups in
$\widetilde{G}\setminus (G_1\setminus\{1_{\times_1}\})$ is a
constant only dependent on $\widetilde{G}\setminus
(G_1\setminus\{1_{\times_1}\})$. Therefore, the length of a
maximal series of normal sub-multi-groups is also a constant, only
dependent on $\widetilde{G}$.

Applying the induction principle, we know that the length of a
maximal series of normal sub-multi-groups of $\widetilde{G}$ is a
constant under an oriented operations
$\overrightarrow{O}(\widetilde{G})$, only dependent on
$\widetilde{G}$ itself.\quad\quad $\natural$

As a special case of Theorem $1.3.7$, we get a consequence in the
following.

\vskip 4mm

\no{\bf Corollary $1.3.4$}(Jordan-H\"{o}lder theorem) \ {\it For a
finite group $G$, the length of its composition series is a
constant, only dependent on $G$.}

\vskip 3mm

Certainly, we can also find other characteristics for multi-groups
similar to group theory, such as those to establish the
decomposition theory for multi-groups similar to the decomposition
theory of abelian groups, to characterize finite generated
multi-groups, $\cdots$, etc.. More observations can be seen in the
finial section of this chapter.

\vskip 6mm

\no{\bf $1.3.3$ Multi-Rings}

\vskip 4mm

\no{\bf Definition $1.3.4$} \ {\it Let
$\widetilde{R}=\bigcup\limits_{i=1}^mR_i$ be a complete
multi-space with a double operation set
$O(\widetilde{R})=\{(+_i,\times_i) , 1\leq i\leq m\}$. If for any
integers $i, j, \ i\not= j, 1\leq i, j\leq m$, $(R_i; +_i,
\times_i)$ is a ring and

$$ (x+_iy)+_jz = x+_i(y+_jz), \  \ \ (x\times_iy)\times_jz = x\times_i(y\times_jz)$$

\no for $\forall x,y,z\in\widetilde{R}$ and

$$x\times_i(y+_jz) = x\times_iy +_jx\times_iz, \  \ \ (y+_jz)\times_ix = y\times_ix +_jz\times_ix$$

\no if all of these operating results exist, then $\widetilde{R}$
is called a multi-ring. If $(R;+_i,\times_i)$ is a field for any
integer $1\leq i\leq m$, then $\widetilde{R}$ is called a
multi-field.}

\vskip 3mm

For a multi-ring $\widetilde{R}=\bigcup\limits_{i=1}^mR_i$, let
$\widetilde{S}\subset\widetilde{R}$ and $O(\widetilde{S})\subset
O(\widetilde{R})$, if $\widetilde{S}$ is also a multi-ring with a
double operation set $O(\widetilde{S})$ , then we call
$\widetilde{S}$ a {\it sub-multi-ring} of $\widetilde{R}$. We get
a criterion for sub-multi-rings in the following.

\vskip 4mm

\no{\bf Theorem $1.3.8$} \ {\it For a multi-ring
$\widetilde{R}=\bigcup\limits_{i=1}^mR_i$, a subset
$\widetilde{S}\subset\widetilde{R}$ with $O(\widetilde{S})\subset
O(\widetilde{R})$ is a sub-multi-ring of $\widetilde{R}$ if and
only if $(\widetilde{S}\bigcap R_k; +_k,\times_k)$ is a subring of
$(R_k; +_k, \times_k)$ or $\widetilde{S}\bigcap R_k=\emptyset$ for
any integer $k, 1\leq k\leq m$.}

\vskip 3mm

{\it Proof} \ For any integer $k, 1\leq k\leq m$, if
$(\widetilde{S}\bigcap R_k; +_k,\times_k)$ is a subring of $(R_k;
+_k,\times_k)$ or $\widetilde{S}\bigcap R_k=\emptyset$, then since
$\widetilde{S}=\bigcup\limits_{i=1}^m(\widetilde{S}\bigcap R_i)$,
we know that $\widetilde{S}$ is a sub-multi-ring by the definition
of a sub-multi-ring.

Now if $\widetilde{S}=\bigcup\limits_{j=1}^sS_{i_j}$ is a
sub-multi-ring of $\widetilde{R}$ with a double operation set
$O(\widetilde{S})=\{(+_{i_j},\times_{i_j}), 1\leq j\leq s\}$, then
$(S_{i_j};  +_{i_j},\times_{i_j})$ is a subring of $(R_{i_j};
+_{i_j}, \times_{i_j})$. Therefore,
$S_{i_j}=R_{i_j}\bigcap\widetilde{S}$ for any integer $j, 1\leq
j\leq s$. But $\widetilde{S}\bigcap S_l=\emptyset$ for other
integer $l\in\{i; 1\leq i\leq m\}\setminus\{i_j; 1\leq j\leq s\}$.
\quad\quad $\natural$

Applying these criterions for subrings of a ring, we get a result
in the following.

\vskip 4mm

\no{\bf Theorem $1.3.9$} \ {\it For a multi-ring
$\widetilde{R}=\bigcup\limits_{i=1}^mR_i$, a subset
$\widetilde{S}\subset\widetilde{R}$ with $O(\widetilde{S})\subset
O(\widetilde{R})$ is a sub-multi-ring of $\widetilde{R}$ if and
only if $(\widetilde{S}\bigcap R_j; +_j) \prec (R_j;+_j)$ and
$(\widetilde{S}; \times_j)$ is complete for any double operation
$(+_j, \times_j)\in O(\widetilde{S})$.}

\vskip 3mm

{\it Proof} \ According to Theorem $1.3.8$, we know that
$\widetilde{S}$ is a sub-multi-ring if and only if
$(\widetilde{S}\bigcap R_i;+_i,\times_i)$ is a subring of
$(R_i;+_i,\times_i)$ or $\widetilde{S}\bigcap R_i=\emptyset$ for
any integer $i, 1\leq i\leq m$. By a well known criterion for
subrings of a ring (see also $[73]$), we know that
$(\widetilde{S}\bigcap R_i;+_i,\times_i)$ is a subring of
$(R_i;+_i,\times_i)$ if and only if $(\widetilde{S}\bigcap R_j;
+_j) \prec (R_j;+_j)$ and $(\widetilde{S}; \times_j)$ is a
complete set for any double operation $(+_j ,\times_j)\in
O(\widetilde{S})$. This completes the proof. \quad\quad $\natural$

We use multi-ideal chains of a multi-ring to characteristic its
structure properties. A {\it multi-ideal} $\widetilde{I}$ of a
multi-ring $\widetilde{R}=\bigcup\limits_{i=1}^mR_i$ with a double
operation set $O(\widetilde{R})$ is a sub-multi-ring of
$\widetilde{R}$ satisfying the following conditions:

\vskip 3mm

$(i)$ \ $\widetilde{I}$ is a sub-multi-group with an operation set
$\{+ | \ (+ ,\times)\in O(\widetilde{I})\}$;

$(ii)$ \ for any $r\in\widetilde{R}, a\in\widetilde{I}$ and
$(+,\times)\in O(\widetilde{I})$, $r\times a\in\widetilde{I}$ and
$a\times r\in\widetilde{I}$ if all of these operating results
exist.

\vskip 4mm

\no{\bf Theorem $1.3.10$} \ {\it A subset $\widetilde{I}$ with
$O(\widetilde{I}), O(\widetilde{I})\subset O(\widetilde{R})$ of a
multi-ring $\widetilde{R}=\bigcup\limits_{i=1}^mR_i$ with a double
operation set $O(\widetilde{R})=\{(+_i,\times_i)| \ 1\leq i\leq
m\}$ is a multi-ideal if and only if $(\widetilde{I}\bigcap R_i,
+_i,\times_i)$ is an ideal of the ring $(R_i,+_i,\times_i)$ or
$\widetilde{I}\bigcap R_i=\emptyset$ for any integer $i, 1\leq
i\leq m$.}

\vskip 3mm

{\it Proof} \ By the definition of a multi-ideal, the necessity of
these conditions is obvious.

For the sufficiency, denote by $\widetilde{R}(+,\times)$ the set
of elements in $\widetilde{R}$ with binary operations ¡°$+$¡± and
¡°$\times$¡±. If there exists an integer $i$ such that
$\widetilde{I}\bigcap R_i\not=\emptyset$ and
$(\widetilde{I}\bigcap R_i, +_i,\times_i)$ is an ideal of
$(R_i,+_i,\times_i)$, then for $\forall a\in\widetilde{I}\bigcap
R_i$, $\forall r_i\in R_i$, we know that

$$r_i\times_i a\in\widetilde{I}\bigcap R_i; \ \ \
a\times_i r_i\in\widetilde{I}\bigcap R_i .$$

Notice that $\widetilde{R}(+_i,\times_i)=R_i$. Thereafter, we get
that

$$r\times_i a\in\widetilde{I}\bigcap R_i \ \ {\rm and} \ \
a\times_i r\in\widetilde{I}\bigcap R_i ,$$

\no for $\forall r\in\widetilde{R}$ if all of these operating
results exist. Whence, $\widetilde{I}$ is a multi-ideal of
$\widetilde{R}$. \quad\quad $\natural$

A multi-ideal $\widetilde{I}$ of a multi-ring $\widetilde{R}$ is
said to be {\it maximal} if for any multi-ideal $\widetilde{I}'$,
$\widetilde{R}\supseteq\widetilde{I}'\supseteq\widetilde{I}$
implies that $\widetilde{I}'=\widetilde{R}$ or
$\widetilde{I}'=\widetilde{I}$. For an order of the double
operations in the set $O(\widetilde{R})$ of a multi-ring
$\widetilde{R}=\bigcup\limits_{i=1}^mR_i$, not loss of generality,
let the order be $(+_1,\times_1)\succ (+_2,\times_2)\succ
\cdots\succ(+_m,\times_m)$, we can define a {\it multi-ideal
chain} of $\widetilde{R}$ by the following programming.

\vskip 3mm

$(i)$ Construct a multi-ideal chain

$$\widetilde{R}\supset\widetilde{R}_{11}\supset\widetilde{R}_{12}
\supset\cdots\supset\widetilde{R}_{1s_1}$$

\no under the double operation $(+_1,\times_1)$, where
$\widetilde{R}_{11}$ is a maximal multi-ideal of $\widetilde{R}$
and in general, $\widetilde{R}_{1(i+1)}$ is a maximal multi-ideal
of $\widetilde{R}_{1i}$ for any integer $i, \ 1\leq i\leq m-1$.

$(ii)$ If a multi-ideal chain

$$\widetilde{R}\supset\widetilde{R}_{11}\supset\widetilde{R}_{12}
\supset\cdots\supset\widetilde{R}_{1s_1}\supset\cdots\supset\widetilde{R}_{i1}
\supset\cdots\supset\widetilde{R}_{is_i}$$

\no has been constructed for $(+_1,\times_1)\succ
(+_2,\times_2)\succ \cdots\succ(+_i,\times_i)$, $1\leq i\leq m-1$,
then construct a multi-ideal chain of $\widetilde{R}_{is_i}$

$$\widetilde{R}_{is_i}\supset\widetilde{R}_{(i+1)1}\supset\widetilde{R}_{(i+1)2}
\supset\cdots\supset\widetilde{R}_{(i+1)s_1}$$

\no under the double operation $(+_{i+1},\times_{i+1})$,  where
$\widetilde{R}_{(i+1)1}$ is a maximal multi-ideal of
$\widetilde{R}_{is_i}$ and in general,
$\widetilde{R}_{(i+1)(i+1)}$ is a maximal multi-ideal of
$\widetilde{R}_{(i+1)j}$ for any integer $j, 1\leq j\leq s_i-1$.
Define a multi-ideal chain of $\widetilde{R}$ under
$(+_1,\times_1)\succ (+_2,\times_2)\succ
\cdots\succ(+_{i+1},\times_{i+1})$ to be

$$\widetilde{R}\supset\widetilde{R}_{11}
\supset\cdots\supset\widetilde{R}_{1s_1}\supset\cdots\supset\widetilde{R}_{i1}
\supset\cdots\supset\widetilde{R}_{is_i}\supset\widetilde{R}_{(i+1)1}
\supset\cdots\supset\widetilde{R}_{(i+1)s_{i+1}}.$$

Similar to multi-groups, we get a result for multi-ideal chains of
a multi-ring in the following.

\vskip 4mm

\no{\bf Theorem $1.3.11$} \ {\it For a multi-ring
$\widetilde{R}=\bigcup\limits_{i=1}^mR_i$, its multi-ideal chain
only has finite terms if and only if the ideal chain of the ring
$(R_i;+_i,\times_i)$ has finite terms, i.e., each ring
$(R_i;+_i,\times_i)$ is an Artin ring for any integer $i, 1\leq
i\leq m$.}

\vskip 3mm

{\it Proof} \ Let the order of these double operations in
$\overrightarrow{O}(\widetilde{R})$ be

$$(+_1,\times_1)\succ (+_2,\times_2)\succ\cdots\succ (+_m,\times_m)$$

\no and let a maximal ideal chain in the ring $(R_1;+_1,\times_1)$
be

$$R_1\succ R_{11}\succ\cdots\succ R_{1t_1}.$$

\no Calculate

$$\widetilde{R}_{11}=\widetilde{R}\setminus\{R_1\setminus R_{11}\}=
R_{11}\bigcup (\bigcup\limits_{i=2}^mR_i),$$

$$\widetilde{R}_{12}=\widetilde{R}_{11}\setminus\{R_{11}\setminus R_{12}\}=
R_{12}\bigcup (\bigcup\limits_{i=2}^mR_i),$$

$$\cdots\cdots\cdots\cdots\cdots\cdots$$

$$\widetilde{R}_{1t_1}=\widetilde{R}_{1t_1}\setminus\{R_{1(t_1-1)}\setminus R_{1t_1}\}=
R_{1t_1}\bigcup (\bigcup\limits_{i=2}^mR_i).$$

\no According to Theorem $1.3.10$, we know that

$$\widetilde{R}\supset\widetilde{R}_{11}\supset\widetilde{R}_{12}
\supset\cdots\supset\widetilde{R}_{1t_1}$$

\no  is a maximal multi-ideal chain of $\widetilde{R}$ under the
double operation $(+_1,\times_1)$. In general, for any integer $i,
1\leq i\leq m-1$, assume

$$R_i\succ R_{i1}\succ\cdots\succ R_{it_i}$$

\no is a maximal ideal chain in the ring
$(R_{(i-1)t_{i-1}};+_i,\times_i)$. Calculate

$$\widetilde{R}_{ik}=R_{ik}\bigcup(\bigcup\limits_{j=i+1}^m\widetilde{R}_{ik}\bigcap R_i).$$

\no Then we know that

$$\widetilde{R}_{(i-1)t_{i-1}}\supset\widetilde{R}_{i1}\supset\widetilde{R}_{i2}
\supset\cdots\supset\widetilde{R}_{it_i}$$

\no is a maximal multi-ideal chain of
$\widetilde{R}_{(i-1)t_{i-1}}$ under the double operation
$(+_i,\times_i)$ by Theorem $3.10$. Whence, if the ideal chain of
the ring $(R_i;+_i,\times_i)$ has finite terms for any integer $i,
1\leq i\leq m$, then the multi-ideal chain of the multi-ring
$\widetilde{R}$ only has finite terms. Now if there exists an
integer $i_0$ such that the ideal chain of the ring
$(R_{i_0},+_{i_0},\times_{i_0})$ has infinite terms, then there
must also be infinite terms in a multi-ideal chain of the
multi-ring $\widetilde{R}$. \quad\quad $\natural$.

A multi-ring is called an {\it Artin multi-ring} if its each
multi-ideal chain only has finite terms. We get a consequence by
Theorem $1.3.11$.

\vskip 4mm

\no{\bf Corollary $1.3.5$} {\it A multi-ring
$\widetilde{R}=\bigcup\limits_{i=1}^mR_i$ with a double operation
set $O(\widetilde{R})=\{(+_i,\times_i)| \ 1\leq i\leq m\}$ is an
Artin multi-ring if and only if the ring $(R_i;+_i,\times_i)$ is
an Artin ring for any integer $i, 1\leq i\leq m$.}

\vskip 3mm

For a multi-ring $\widetilde{R}=\bigcup\limits_{i=1}^mR_i$ with a
double operation set $O(\widetilde{R})=\{(+_i,\times_i)| \ 1\leq
i\leq m\}$, an element $e$ is an {\it idempotent} element if
$e_{\times}^2 = e\times e = e$ for a double binary operation
$(+,\times)\in O(\widetilde{R})$. We define the {\it directed sum}
$\widetilde{I}$ of two multi-ideals $\widetilde{I}_1$ and
$\widetilde{I}_2$ as follows:

\vskip 3mm

$(i)$ $\widetilde{I}=\widetilde{I}_1\bigcup\widetilde{I}_2$;

$(ii)$ $\widetilde{I}_1\bigcap\widetilde{I}_2=\{0_+\}, \ {\rm or}
\ \widetilde{I}_1\bigcap\widetilde{I}_2=\emptyset$, where $0_+$
denotes an unit element under the operation $+$.

Denote the directed sum of $\widetilde{I}_1$ and $\widetilde{I}_2$
by

$$\widetilde{I}=\widetilde{I}_1\bigoplus\widetilde{I}_2.$$

If $\widetilde{I}=\widetilde{I}_1\bigoplus\widetilde{I}_2$ for any
$\widetilde{I}_1, \widetilde{I}_2$ implies that
$\widetilde{I}_1=\widetilde{I}$ or
$\widetilde{I}_2=\widetilde{I}$, then $\widetilde{I}$ is called
{\it non-reducible}. We get the following result for Artin
multi-rings similar to a well-known result for Artin rings (see
$[107]$ for details).

\vskip 4mm

\no{\bf Theorem $1.3.12$} \ {\it Any Artin multi-ring
$\widetilde{R}=\bigcup\limits_{i=1}^mR_i$ with a double operation
set $O(\widetilde{R})=\{(+_i,\times_i)| \ 1\leq i\leq m\}$ is a
directed sum of finite non-reducible multi-ideals, and if
$(R_i;+_i,\times_i)$ has unit $1_{\times_i}$ for any integer
$i,1\leq i\leq m$, then}

$$\widetilde{R}=\bigoplus\limits_{i=1}^m(\bigoplus\limits_{j=1}^{s_i}(R_i\times_ie_{ij})
\bigcup (e_{ij}\times_i R_i)),$$

\no{\it where $e_{ij}, 1\leq j\leq s_i$ are orthogonal idempotent
elements of the ring $R_i$.}

\vskip 3mm

{\it Proof} \ Denote by $\widetilde{M}$ the set of multi-ideals
which can not be represented by a directed sum of finite
multi-ideals in $\widetilde{R}$. According to Theorem $3.11$,
there is a minimal multi-ideal $\widetilde{I}_0$ in
$\widetilde{M}$. It is obvious that $\widetilde{I}_0$ is
reducible.

Assume that $\widetilde{I}_0=\widetilde{I}_1+\widetilde{I}_2$.
Then $\widetilde{I}_1\not\in\widetilde{M}$ and
$\widetilde{I}_2\not\in\widetilde{M}$. Therefore,
$\widetilde{I}_1$ and $\widetilde{I}_2$ can be represented by a
directed sum of finite multi-ideals. Thereby $\widetilde{I}_0$ can
be also represented  by a directed sum of finite multi-ideals.
Contradicts that $\widetilde{I}_0\in\widetilde{M}$.

Now let

$$\widetilde{R}=\bigoplus\limits_{i=1}^s\widetilde{I}_i,$$

\no where each $\widetilde{I}_i, 1\leq i\leq s$ is non-reducible.
Notice that for a double operation $(+,\times )$, each
non-reducible multi-ideal of $\widetilde{R}$ has the form

$$(e\times R(\times))\bigcup (R(\times )\times e), \ \ e\in R(\times ).$$

\no Whence, we know that there is a set $T\subset\widetilde{R}$
such that

$$\widetilde{R}=\bigoplus\limits_{e\in T, \ \times\in O(\widetilde{R})}
(e\times R(\times))\bigcup (R(\times )\times e).$$

For any operation $\times\in O(\widetilde{R})$ and the unit
$1_{\times}$, assume that

$$1_{\times} = e_1\oplus e_2\oplus\cdots \oplus e_l, \ e_i\in T, \ 1\leq i\leq s.$$

\no Then

$$e_i\times 1_{\times}= (e_i\times e_1)\oplus (e_i\times e_2)\oplus\cdots \oplus (e_i\times e_l).$$

\no Therefore, we get that

$$e_i = e_i\times e_i=e_i^2 \ \ {\rm and} \ \ e_i\times e_j=0_i \ \ {\rm for} \ \ i\not=j.$$

\no That is, $e_i, 1\leq i\leq l$ are orthogonal idempotent
elements of $\widetilde{R}(\times)$. Notice that
$\widetilde{R}(\times)=R_h$ for some integer $h$. We know that
$e_i, 1\leq i\leq l$ are orthogonal idempotent elements of the
ring $(R_h,+_h,\times_h)$. Denote by $e_{hi}$ for $e_i$, $1\leq
i\leq l$. Consider all units in $\widetilde{R}$, we get that

$$\widetilde{R}=\bigoplus\limits_{i=1}^m(\bigoplus\limits_{j=1}^{s_i}(R_i\times_ie_{ij})
\bigcup (e_{ij}\times_i R_i)).$$

\no This completes the proof. \quad\quad $\natural$

\vskip 4mm

\no{\bf Corollary $1.3.6$} \ {\it Any Artin ring $(R \ ;+,\times
)$ is a directed sum of finite ideals, and if $(R \ ;+,\times )$
has unit $1_{\times}$, then}

$$R =\bigoplus\limits_{i=1}^s R_ie_{i },$$

\no{\it where $e_{i}, 1\leq i\leq s$ are orthogonal idempotent
elements of the ring $(R;+,\times )$.}

\vskip 3mm

Similarly, we can also define {\it Noether multi-rings, simple
multi-rings, half-simple multi-rings}, $\cdots$, etc. and find
their algebraic structures.

\vskip 6mm

\no{\bf $1.3.4$ \ Multi-Vector spaces}

\vskip 4mm

\no{\bf Definition $1.3.5$} \ {\it Let
$\widetilde{V}=\bigcup\limits_{i=1}^k V_i$ be a complete
multi-space with an operation set $O(\widetilde{V})=\{
(\dot{+}_i,\cdot_i) \ | \ 1\leq i\leq m\}$ and let
$\widetilde{F}=\bigcup\limits_{i=1}^k F_i$ be a multi-filed with a
double operation set $O(\widetilde{F})=\{(+_i,\times_i )\ | \
1\leq i\leq k\}$. If for any integers $i,j, \ 1\leq i, j\leq k$
and $\forall {\bf a,b,c}\in\widetilde{V}$,
$k_1,k_2\in\widetilde{F}$,

$(i)$ $(V_i;\dot{+}_i,\cdot_i)$ is a vector space on $F_i$ with
vector additive ¡°$\dot{+}_i$¡± and scalar multiplication
¡°$\cdot_i$¡±;

$(ii)$ $({\bf a}\dot{+}_i{\bf b})\dot{+}_j{\bf c}= {\bf
a}\dot{+}_i({\bf b}\dot{+}_j{\bf c})$;

$(iii)$ $(k_1+_i k_2)\cdot_j{\bf a}=k_1+_i(k_2\cdot_j{\bf a});$

\no provided these operating results exist, then $\widetilde{V}$
is called a multi-vector space on the multi-filed space
$\widetilde{F}$ with an double operation set $O(\widetilde{V})$,
denoted by $(\widetilde{V}; \widetilde{F})$.}

\vskip 3mm

For subsets $\widetilde{V}_1\subset\widetilde{V}$ and
$\widetilde{F}_1\subset\widetilde{F}$, if $(\widetilde{V}_1;
\widetilde{F}_1)$ is also a multi-vector space, then we call
$(\widetilde{V}_1; \widetilde{F}_1)$ a {\it multi-vector subspace}
of $(\widetilde{V}; \widetilde{F})$. Similar to the linear space
theory, we get the following criterion for multi-vector subspaces.

\vskip 4mm

\no{\bf Theorem $1.3.13$} \ {\it For a multi-vector space
$(\widetilde{V}; \widetilde{F})$,
$\widetilde{V}_1\subset\widetilde{V}$ and
$\widetilde{F}_1\subset\widetilde{F}$, $(\widetilde{V}_1;
\widetilde{F}_1)$ is a multi-vector subspace of $(\widetilde{V};
\widetilde{F})$ if and only if for any vector additive
¡°$\dot{+}$¡±, scalar multiplication ¡°$\cdot$¡± in
$(\widetilde{V}_1; \widetilde{F}_1)$ and $\forall {\bf
a,b}\in\widetilde{V}$, $\forall \alpha\in\widetilde{F}$,}

$$\alpha\cdot {\bf a}\dot{+}{\bf b}\in\widetilde{V}_1$$

\no{\it provided these operating results exist.}

\vskip 3mm

{\it Proof} \ Denote by $\widetilde{V}=\bigcup\limits_{i=1}^k V_i,
\widetilde{F}=\bigcup\limits_{i=1}^k F_i$. Notice that
$\widetilde{V}_1=\bigcup\limits_{i=1}^k(\widetilde{V}_1\bigcap
V_i)$. By definition, we know that  $(\widetilde{V}_1;
\widetilde{F}_1)$ is a multi-vector subspace of $(\widetilde{V};
\widetilde{F})$ if and only if for any integer $i, 1\leq i\leq k$,
$(\widetilde{V}_1\bigcap V_i; \dot{+}_i,\cdot_i)$ is a vector
subspace of $(V_i, \dot{+}_i,\cdot_i)$ and $\widetilde{F}_1$ is a
multi-filed subspace of $\widetilde{F}$ or $\widetilde{V}_1\bigcap
V_i=\emptyset$.

According to a criterion for linear subspaces of a linear space
([$33$]), we know that $(\widetilde{V}_1\bigcap V_i;
\dot{+}_i,\cdot_i)$ is a vector subspace of $(V_i,
\dot{+}_i,\cdot_i)$ for any integer $i, 1\leq i\leq k$ if and only
if for $\forall {\bf a, b}\in\widetilde{V}_1\bigcap V_i$,
$\alpha\in F_i$,

$$\alpha\cdot_i{\bf a}\dot{+}_i{\bf b}\in\widetilde{V}_1\bigcap V_i.$$

\no That is, for any vector additive ¡°$\dot{+}$¡±, scalar
multiplication¡°$\cdot$¡±in $(\widetilde{V}_1; \widetilde{F}_1)$
and $\forall {\bf a,b}\in\widetilde{V}$, $\forall
\alpha\in\widetilde{F}$, if $\alpha\cdot {\bf a}\dot{+}{\bf b}$
exists, then $\alpha\cdot {\bf a}\dot{+}{\bf
b}\in\widetilde{V}_1$.\quad\quad $\natural$

\vskip 4mm

\no{\bf Corollary $1.3.7$} \ {\it Let
$(\widetilde{U};\widetilde{F}_1),(\widetilde{W};\widetilde{F}_2)$
be two multi-vector subspaces of a multi-vector space
$(\widetilde{V};\widetilde{F})$. Then
$(\widetilde{U}\bigcap\widetilde{W};\widetilde{F}_1\bigcap\widetilde{F}_2)$
is a multi-vector space.}

\vskip 3mm

For a multi-vector space $(\widetilde{V};\widetilde{F})$, vectors
${\bf a}_1,{\bf a}_2,\cdots ,{\bf a}_n\in\widetilde{V}$, if there
are scalars $\alpha_1, \alpha_2, \cdots ,\alpha_n\in\widetilde{F}$
such that

$$\alpha_1\cdot_1{\bf a}_1\dot{+}_1\alpha_2\cdot_2{\bf a}_2\dot{+}_2\cdots\dot{+}_{n-1}\alpha_n\cdot_n
{\bf a}_n={\bf 0},$$

\no where ${\bf 0}\in\widetilde{V}$ is a unit under an operation
¡°$+$¡± in $\widetilde{V}$ and $\dot{+}_i, \cdot_i\in
O(\widetilde{V})$, then these vectors ${\bf a}_1,{\bf a}_2,\cdots
,{\bf a}_n$ are said to be {\it linearly dependent}. Otherwise,
${\bf a}_1,{\bf a}_2,\cdots ,{\bf a}_n$ are said to be {\it
linearly independent}.

Notice that there are two cases for linearly independent vectors
${\bf a}_1,{\bf a}_2,\cdots ,{\bf a}_n$ in a multi-vector space:

\vskip 2mm

($i$) \ for scalars $\alpha_1, \alpha_2, \cdots
,\alpha_n\in\widetilde{F}$, if

$$\alpha_1\cdot_1{\bf a}_1\dot{+}_1\alpha_2\cdot_2{\bf a}_2
\dot{+}_2\cdots\dot{+}_{n-1}\alpha_n\cdot_n {\bf a}_n={\bf 0},$$

\no where ${\bf 0}$ is a unit of $\widetilde{V}$ under an
operation ¡°$+$¡± in $O(\widetilde{V})$, then $\alpha_1=0_{+_1},
\alpha_2=0_{+_2}, \cdots , \alpha_n=0_{+_n}$, where $0_{+_i}$
 is the unit under the operation ¡°$+_i$¡±in
$\widetilde{F}$ for integer $i, 1\leq i\leq n$.

($ii$) \ the operating result of $\alpha_1\cdot_1{\bf
a}_1\dot{+}_1\alpha_2\cdot_2{\bf
a}_2\dot{+}_2\cdots\dot{+}_{n-1}\alpha_n\cdot_n {\bf a}_n$ does
not exist.\vskip 2mm

Now for a subset $\widehat{S}\subset\widetilde{V}$, define its
{\it linearly spanning set} $\left<\widehat{S}\right>$ to be

$$\left<\widehat{S}\right> = \{ \ {\bf a} \ | \ {\bf a}=\alpha_1\cdot_1{\bf a}_1\dot{+}_1
\alpha_2\cdot_2{\bf a}_2\dot{+}_2\cdots \in\widetilde{V}, {\bf
a}_i\in\widehat{S}, \alpha_i\in\widetilde{F}, i\geq 1 \}.$$

\no For a multi-vector space $(\widetilde{V};\widetilde{F})$, if
there exists a subset $\widehat{S},
\widehat{S}\subset\widetilde{V}$ such that $\widetilde{V} =
\left<\widehat{S}\right>$, then we say $\widehat{S}$ is a {\it
linearly spanning set} of the multi-vector space $\widetilde{V}$.
If these vectors in a linearly spanning set $\widehat{S}$ of the
multi-vector space $\widetilde{V}$ are linearly independent, then
$\widehat{S}$ is said to be a {\it basis} of $\widetilde{V}$.

\vskip 4mm

\no{\bf Theorem $1.3.14$} \ {Any multi-vector space
$(\widetilde{V};\widetilde{F})$ has a basis.}

\vskip 3mm

{\it Proof} \ Assume $\widetilde{V}=\bigcup\limits_{i=1}^k V_i,
\widetilde{F}=\bigcup\limits_{i=1}^k F_i$ and the basis of the
vector space $(V_i;\dot{+}_i,\cdot_i)$ is $\Delta_i=\{{\bf
a}_{i1},{\bf a}_{i2},\cdots ,{\bf a}_{in_i}\}$, $1\leq i\leq k$.
Define

$$\widehat{\Delta} \ = \ \bigcup\limits_{i=1}^k\Delta_i.$$

\no Then $\widehat{\Delta}$ is a linearly spanning set for
$\widetilde{V}$ by definition.

If these vectors in $\widehat{\Delta}$ are linearly independent,
then $\widehat{\Delta}$ is a basis of $\widetilde{V}$. Otherwise,
choose a vector ${\bf b}_1\in\widehat{\Delta}$ and define
$\widehat{\Delta}_1=\widehat{\Delta}\setminus\{{\bf b}_1\}$.

If we have obtained a set $\widehat{\Delta}_s, s\geq 1$ and it is
not a basis, choose a vector ${\bf b}_{s+1}\in\widehat{\Delta}_s$
and define
$\widehat{\Delta}_{s+1}=\widehat{\Delta}_s\setminus\{{\bf
b}_{s+1}\}$.

If these vectors in $\widehat{\Delta}_{s+1}$ are linearly
independent, then $\widehat{\Delta}_{s+1}$ is a basis of
$\widetilde{V}$. Otherwise, we can define a set
$\widehat{\Delta}_{s+2}$ again. Continue this process. Notice that
all vectors in $\Delta_i$ are linearly independent for any integer
$i, 1\leq i\leq k$. Therefore, we can finally get a basis of
$\widetilde{V}$. \quad\quad $\natural$

Now we consider finite-dimensional multi-vector spaces. A
multi-vector space $\widetilde{V}$ is {\it finite-dimensional} if
it has a finite basis. By Theorem $1.2.14$, if the vector space
$(V_i; +_i,\cdot_i)$ is finite-dimensional for any integer $i,
1\leq i\leq k$, then $(\widetilde{V};\widetilde{F})$ is
finite-dimensional. On the other hand, if there is an integer
$i_0, 1\leq i_0\leq k$ such that the vector space $(V_{i_0};
+_{i_0},\cdot_{i_0})$ is infinite-dimensional, then
$(\widetilde{V};\widetilde{F})$ is also infinite-dimensional. This
enables us to get a consequence in the following.

\vskip 4mm

\no{\bf Corollary $1.3.8$} \ {\it Let
$(\widetilde{V};\widetilde{F})$ be a multi-vector space with
$\widetilde{V}=\bigcup\limits_{i=1}^k V_i,
\widetilde{F}=\bigcup\limits_{i=1}^k F_i$. Then
$(\widetilde{V};\widetilde{F})$ is finite-dimensional if and only
if $(V_i; +_i,\cdot_i)$ is finite-dimensional for any integer $i,
1\leq i\leq k$.}

\vskip 4mm

\no{\bf Theorem $1.3.15$} \ {\it For a finite-dimensional
multi-vector space $(\widetilde{V};\widetilde{F})$, any two bases
have the same number of vectors.}

\vskip 3mm

{\it Proof} \ Let $\widetilde{V}=\bigcup\limits_{i=1}^k V_i$ and
$\widetilde{F}=\bigcup\limits_{i=1}^k F_i$. The proof is by the
induction on $k$. For $k=1$, the assertion is true by Theorem $4$
of Chapter $2$ in $[33]$.

For the case of $k=2$, notice that by a result in linearly vector
spaces (see also [$33$]), for two subspaces $W_1,W_2$ of a
finite-dimensional vector space, if the basis of $W_1\bigcap W_2$
is $\{{\bf a}_1,{\bf a}_2,\cdots , {\bf a}_t\}$, then the basis of
$W_1\bigcup W_2$ is

$$\{{\bf a}_1,{\bf a}_2,\cdots , {\bf a}_t, {\bf b}_{t+1},
{\bf b}_{t+2},\cdots ,{\bf b}_{dimW_1}, {\bf c}_{t+1},{\bf
c}_{t+2}, \cdots ,{\bf c}_{dimW_2}\},$$

\no where, $\{{\bf a}_1,{\bf a}_2,\cdots , {\bf a}_t, {\bf
b}_{t+1},{\bf b}_{t+2},\cdots ,{\bf b}_{dimW_1}\}$ is a basis of
$W_1$ and $\{{\bf a}_1,{\bf a}_2,\cdots , {\bf a}_t,$ ${\bf
c}_{t+1},{\bf c}_{t+2},\cdots ,{\bf c}_{dimW_2}\}$ a basis of
$W_2$.

Whence, if $\widetilde{V} = W_1\bigcup W_2$ and
$\widetilde{F}=F_1\bigcup F_2$, then the basis of $\widetilde{V}$
is also

$$\{{\bf a}_1,{\bf a}_2,\cdots , {\bf a}_t, {\bf b}_{t+1},
{\bf b}_{t+2},\cdots ,{\bf b}_{dimW_1}, {\bf c}_{t+1},{\bf
c}_{t+2}, \cdots ,{\bf c}_{dimW_2}\}.$$

Assume the assertion is true for $k=l, l\geq 2$. Now we consider
the case of $k=l+1$. In this case, since

$$\widetilde{V}=(\bigcup\limits_{i=1}^l V_i)\bigcup V_{l+1}, \ \widetilde{F}
=(\bigcup\limits_{i=1}^l F_i)\bigcup F_{l+1},$$

\no by the induction assumption, we know that any two bases of the
multi-vector space $(\bigcup\limits_{i=1}^l
V_i;\bigcup\limits_{i=1}^l F_i)$ have the same number $p$ of
vectors. If the basis of $(\bigcup\limits_{i=1}^l V_i)\bigcap
V_{l+1}$ is $\{{\bf e}_1,{\bf e}_2,\cdots ,{\bf e}_n\}$, then the
basis of $\widetilde{V}$ is

$$\{{\bf e}_1,{\bf e}_2,\cdots ,{\bf e}_n, {\bf f}_{n+1},
{\bf f}_{n+2},\cdots ,{\bf f}_{p}, {\bf g}_{n+1},{\bf
g}_{n+2},\cdots ,{\bf g}_{dimV_{l+1}}\},$$

\no where $\{{\bf e}_1,{\bf e}_2,\cdots ,{\bf e}_n, {\bf
f}_{n+1},{\bf f}_{n+2},\cdots ,{\bf f}_{p}\}$ is a basis of
$(\bigcup\limits_{i=1}^l V_i; \bigcup\limits_{i=1}^l F_i)$ and
$\{{\bf e}_1,{\bf e}_2,\cdots ,{\bf e}_n,$ ${\bf g}_{n+1},{\bf
g}_{n+2},\cdots ,{\bf g}_{dimV_{l+1}}\}$ is a basis of $V_{l+1}$.
Whence, the number of vectors in a basis of $\widetilde{V}$ is
$p+dimV_{l+1}-n$ for the case $n=l+1$.

Therefore, we know the assertion is true for any integer $k$ by
the induction principle.\quad\quad $\natural$

The cardinal number of a basis of a finite dimensional
multi-vector space $\widetilde{V}$ is called its {\it dimension},
denoted by $dim\widetilde{V}$.

\vskip 4mm

\no{\bf Theorem $1.3.16$}({\it dimensional formula}) \ {\it For a
multi-vector space $(\widetilde{V}; \widetilde{F})$ with
$\widetilde{V}=\bigcup\limits_{i=1}^k V_i$ and
$\widetilde{F}=\bigcup\limits_{i=1}^k F_i$, the dimension
$dim\widetilde{V}$ of $\widetilde{V}$ is}

$$dim\widetilde{V}=\sum\limits_{i=1}^k(-1)^{i-1}
\sum\limits_{\{i1,i2,\cdots ,ii\}\subset\{1,2,\cdots
,k\}}dim(V_{i1} \bigcap V_{i2}\bigcap\cdots\bigcap V_{ii}).$$

\vskip 3mm

{\it Proof} \ The proof is by the induction on $k$. For $k=1$, the
formula is turn to a trivial case of $dim\widetilde{V}=dimV_1$.
for $k=2$, the formula is

$$dim\widetilde{V}=dimV_1+dimV_2-dim(V_1\bigcap dimV_2),$$

\no which is true by the proof of Theorem $1.3.15$.

Now we assume the formula is true for $k=n$. Consider the case of
$k=n+1$. According to the proof of Theorem $1.3.15$, we know that

\begin{eqnarray*}
dim\widetilde{V} &=& dim(\bigcup\limits_{i=1}^nV_i)+dimV_{n+1}
-dim((\bigcup\limits_{i=1}^nV_i)\bigcap V_{n+1})\\
 &=& dim(\bigcup\limits_{i=1}^nV_i)+dimV_{n+1}
-dim(\bigcup\limits_{i=1}^n(V_i\bigcap V_{n+1}))\\
&=& dimV_{n+1}+
\sum\limits_{i=1}^n(-1)^{i-1}\sum\limits_{\{i1,i2,\cdots,ii\}\subset\{1,2,\cdots
,n\}}dim(V_{i1}\bigcap V_{i2}\bigcap\cdots\bigcap V_{ii})\\
&+&
\sum\limits_{i=1}^n(-1)^{i-1}\sum\limits_{\{i1,i2,\cdots,ii\}\subset\{1,2,\cdots
,n\}}dim(V_{i1}\bigcap V_{i2}\bigcap\cdots\bigcap V_{ii}\bigcap
V_{n+1})\\
&=& \sum\limits_{i=1}^n(-1)^{i-1} \sum\limits_{\{i1,i2,\cdots
,ii\}\subset\{1,2,\cdots ,k\}}dim(V_{i1} \bigcap
V_{i2}\bigcap\cdots\bigcap V_{ii}).
\end{eqnarray*}

\no By the induction principle, we know the formula is true for
any integer $k$. \quad\quad $\natural$

As a consequence, we get the following formula.

\vskip 4mm

\no{\bf Corollary $1.3.9$}({\it additive formula}) \ {\it For any
two multi-vector spaces $\widetilde{V}_1,\widetilde{V}_2$,}

$$dim(\widetilde{V}_1\bigcup\widetilde{V}_2)=dim\widetilde{V}_1+
dim\widetilde{V}_2-dim(\widetilde{V}_1\bigcap\widetilde{V}_2).$$

\vskip 5mm

\no{\bf \S $1.4$ \ Multi-Metric Spaces}

\vskip 3mm

\no{\bf $1.4.1.$ Metric spaces}

\vskip 3mm

\no A set $M$ associated with a metric function $\rho: M\times
M\rightarrow R^+=\{x \ | \ x\in R, x\geq 0\}$ is called a {\it
metric space} if for $\forall x,y,z\in M$, the following
conditions for $\rho$ hold:

\vskip 3mm

($1$)({\it definiteness}) \  $\rho (x,y)=0$ if and only if $x=y$;

($ii$)({\it symmetry}) \ $\rho (x,y)=\rho (y,x)$;

($iii$)({\it triangle inequality}) \ $\rho (x,y)+\rho (y,z)\geq
\rho (x,z).$

\vskip 2mm

A metric space $M$ with a metric function $\rho$ is usually
denoted by $(M \ ;\rho)$. Any $x, x\in M$ is called a point of $(M
\ ;\rho)$. A sequence $\{x_n\}$ is said to be {\it convergent to
$x$} if for any number $\epsilon > 0$ there is an integer $N$ such
that $n\geq N$ implies $\rho (x_n,x) \ < \ 0$, denoted by
$\lim\limits_{n}x_n=x$. We have known the following result in
metric spaces.

\vskip 4mm

\no{\bf Theorem $1.4.1$}\ {\it Any sequence $\{x_n\}$ in a metric
space has at most one limit point.}

\vskip 3mm

For $x_0\in M$ and $\epsilon > 0$, a $\epsilon$-disk about $x_0$
is defined by

$$B(x_0, \epsilon)=\{ \ x \ | \ x\in M, \rho (x,x_0) \ < \ \epsilon \}.$$

\no If $A\subset M$ and there is an $\epsilon$-disk
$B(x_0,\epsilon )\supset A$, we say $A$ is a bounded point set of
$M$.

\vskip 4mm

\no{\bf Theorem $1.4.2$} \ {\it Any convergent sequence $\{x_n\}$
in a metric space is a bounded point set.}

\vskip 3mm

Now let $(M,\rho )$ be a metric space and $\{x_n\}$ a sequence in
$M$. If for any number $\epsilon > 0, \epsilon\in {\bf R}$, there
is an integer $N$ such that $n,m\geq N$ implies $\rho (x_n,x_m) \
< \ \epsilon$, we call $\{x_n\}$ a {\it Cauchy sequence}. A metric
space $(M,\rho )$ is called to be {\it completed} if its every
Cauchy sequence converges.

\vskip 4mm

\no{\bf Theorem $1.4.3$} \ {\it For a completed metric space
$(M,\rho )$, if an $\epsilon$-disk sequence $\{B_n\}$ satisfies}

($i$) \ {\it $B_1\supset B_2\supset\cdots\supset
B_n\supset\cdots$};

($ii$) \ $\lim\limits_{n}\epsilon_n=0$,

\no {\it where $\epsilon_n
> 0$ and $B_n=\{ \ x \ | \ x\in M, \rho (x,
x_n)\leq\epsilon_n\}$ for any integer $n, n=1,2,\cdots $, then
$\bigcap\limits_{n=1}^{\infty}B_n$ only has one point.}

\vskip 3mm

For a metric space $(M,\rho )$ and $T: M\rightarrow M$ a mapping
on $(M,\rho )$, if there exists a point $x^*\in M$ such that

$$Tx^*=x^*,$$

\no then $x^*$ is called a {\it fixed point} of $T$. If there
exists a constant $\eta , 0 < \eta  <  1$ such that

$$\rho (Tx, Ty)\leq \eta\rho (x,y)$$

\no for $\forall x,y\in M$, then $T$ is called a {\it
contraction}.

\vskip 4mm

\no{\bf Theorem $1.4.4$} (Banach) \ {\it Let $(M,\rho )$ be a
completed metric space and let $T: M\rightarrow M$ be a
contraction. Then $T$ has only one fixed point.}

\vskip 5mm

\no{\bf $1.4.2.$ Multi-Metric spaces}

\vskip 4mm

\no{\bf Definition $1.4.1$} \ {\it A multi-metric space is a union
$\widetilde{M}=\bigcup\limits_{i=1}^m M_i$ such that each $M_i$ is
a space with a metric $\rho_i$ for $\forall i, 1\leq i\leq m$.}

When we say a {\it multi-metric space
$\widetilde{M}=\bigcup\limits_{i=1}^m M_i$}, it means that a
multi-metric space with metrics $\rho_1,\rho_2, \cdots ,\rho_m$
such that $(M_i,\rho_i)$ is a metric space for any integer $i,
1\leq i\leq m$. For a multi-metric space
$\widetilde{M}=\bigcup\limits_{i=1}^m M_i$, $x\in\widetilde{M}$
and a positive number $R$, a {\it $R$-disk} $B(x,R)$ in
$\widetilde{M}$ is defined by

$$B(x,R)=\{ \ y \ | \ {\rm there \ exists \ an \
integer} \ k, 1\leq k\leq m \ {\rm such \ that} \ \rho_k(y,x) \ <
R, y\in\widetilde{M}\}$$

\vskip 3mm

\no{\bf Remark $1.4.1$} \ The following two extremal cases are
permitted in Definition $1.4.1$:

($i$) \ there are integers $i_1,i_2,\cdots,i_s$ such that
$M_{i_1}=M_{i_2}=\cdots =M_{i_s}$, where $i_j\in\{1,2,\cdots ,
m\}$, $1\leq j\leq s$;

($ii$) \ there are integers $l_1,l_2,\cdots,l_s$ such that
$\rho_{l_1}=\rho_{l_2}=\cdots =\rho_{l_s}$, where
$l_j\in\{1,2,\cdots , m\}$, $1\leq j\leq s$.

For metrics on a space, we have the following result.

\vskip 4mm

\no{\bf Theorem $1.4.5$} \ {\it Let $\rho_1,\rho_2,\cdots ,\rho_m$
be $m$ metrics on a space $M$ and let $F$ be a function on ${\bf
R}^m$ such that the following conditions hold:}

($i$) \ {\it $F(x_1,x_2,\cdots ,x_m)\geq F(y_1,y_2,\cdots ,y_m)$
for $\forall i, 1\leq i\leq m$, $x_i\geq y_i$};

($ii$) \ {\it $F(x_1,x_2,\cdots ,x_m)=0$ only if $x_1=x_2=\cdots
=x_m=0$};

($iii$) \ {\it for two $m$-tuples $(x_1,x_2,\cdots ,x_m)$ and
$(y_1,y_2,\cdots ,y_m)$},

$$F(x_1,x_2,\cdots ,x_m)+F(y_1,y_2,\cdots ,y_m)\geq
F(x_1+y_1,x_2+y_2,\cdots ,x_m+y_m).$$

\no {\it Then $F(\rho_1,\rho_2,\cdots ,\rho_m)$ is also a metric
on $M$. }

\vskip 3mm

{\it Proof} \ We only need to prove that $F(\rho_1,\rho_2,\cdots
,\rho_m)$ satisfies those of metric conditions for $\forall
x,y,z\in M$.

By ($ii$), $F(\rho_1(x,y),\rho_2(x,y),\cdots ,\rho_m(x,y))=0$ only
if $\rho_i(x,y)=0$ for any integer $i$. Since $\rho_i$ is a metric
on $M$, we know that $x=y$.

For any integer $i, 1\leq i\leq m$, since $\rho_i$ is a metric on
$M$, we know that $\rho_i(x,y)=\rho_i(y,x)$. Whence,

$$F(\rho_1(x,y),\rho_2(x,y),\cdots ,\rho_m(x,y))=
F(\rho_1(y,x),\rho_2(y,x),\cdots ,\rho_m(y,x)).$$

Now by ($i$) and ($iii$), we get that

\begin{eqnarray*}
& \ & F(\rho_1(x,y),\rho_2(x,y),\cdots
,\rho_m(x,y))+F(\rho_1(y,z),\rho_2(y,z),\cdots ,\rho_m(y,z))\\
& \ &\geq
F(\rho_1(x,y)+\rho_1(y,z),\rho_2(x,y)+\rho_2(y,z),\cdots,\rho_m(x,y)+\rho_m(y,z))\\
& \ &\geq F(\rho_1(x,z),\rho_2(x,z),\cdots ,\rho_m(x,z)).
\end{eqnarray*}

\no Therefore, $F(\rho_1,\rho_2,\cdots ,\rho_m)$ is a metric on
$M$. \quad\quad $\natural$

\vskip 4mm

\no{\bf Corollary $1.4.1$} \ {If $\rho_1,\rho_2,\cdots ,\rho_m$
are $m$ metrics on a space $M$, then $\rho_1+\rho_2+\cdots
+\rho_m$ and
$\frac{\rho_1}{1+\rho_1}+\frac{\rho_2}{1+\rho_2}+\cdots
+\frac{\rho_m}{1+\rho_m}$ are also metrics on $M$.}

\vskip 3mm

A sequence $\{x_n\}$ in a multi-metric space
$\widetilde{M}=\bigcup\limits_{i=1}^m M_i$ is said to be {\it
convergent to a point} $x, x\in\widetilde{M}$ if for any number
$\epsilon > 0$, there exist numbers $N$ and $i, 1\leq i\leq m$
such that

$$\rho_i(x_n,x) \ < \ \epsilon$$

\no provided $n\geq N$. If $\{x_n\}$ is convergent to a point $x,
x\in\widetilde{M}$, we denote it by $\lim\limits_{n}x_n = x$.

We get a characteristic for convergent sequences in a multi-metric
space as in the following.

\vskip 4mm

\no{\bf Theorem $1.4.6$} \ {\it A sequence $\{x_n\}$ in a
multi-metric space $\widetilde{M}=\bigcup\limits_{i=1}^m M_i$ is
convergent if and only if there exist integers $N$ and $k, 1\leq
k\leq m$ such that the subsequence $\{x_n| n\geq N\}$ is a
convergent sequence in $(M_k,\rho_k)$.}

\vskip 3mm

{\it Proof} \ If there exist integers $N$ and $k,1\leq k\leq m$
such that $\{x_n | n\geq N\}$ is a convergent sequence in
$(M_k,\rho_k)$, then for any number $\epsilon  > 0$, by definition
there exist an integer $P$ and a point $x, x\in M_k$ such that

$$\rho_k(x_n,x) \ < \ \epsilon$$

\no if $n\geq max\{ N, \ P \}$.

Now if $\{x_n\}$ is a convergent sequence in the multi-space
$\widetilde{M}$, by definition for any positive number $\epsilon
> 0$, there exist a point $x, x\in\widetilde{M}$, natural numbers
$N(\epsilon )$ and integer $k, 1\leq k\leq m$ such that if $n\geq
N(\epsilon)$, then

$$\rho_k(x_n,x) \ < \ \epsilon.$$

\no That is, $\{ x_n | n\geq N(\epsilon) \}\subset M_k$ and $\{
x_n | n\geq N(\epsilon) \}$ is a convergent sequence in
$(M_k,\rho_k)$.\quad\quad $\natural$

\vskip 4mm

\no{\bf Theorem $1.4.7$} \ {\it Let
$\widetilde{M}=\bigcup\limits_{i=1}^m M_i$ be a multi-metric
space. For two sequences $\{x_n\},$ $ \{y_n\}$ in $\widetilde{M}$,
if $\lim\limits_{n}x_n = x_0$, $\lim\limits_{n}y_n = y_0$ and
there is an integer $p$ such that $x_0,y_0\in M_p$, then
$\lim\limits_n\rho_p(x_n,y_n)=\rho_p(x_0,y_0)$.}

\vskip 3mm

{\it Proof} \ According to Theorem $1.4.6$, there exist integers
$N_1$ and $N_2$ such that if $n\geq max\{N_1, N_2\}$, then
$x_n,y_n\in M_p$. Whence, we know that

$$\rho_p(x_n,y_n) \leq \rho_p(x_n,x_0)+\rho_p(x_0,y_0)+\rho_p(y_n,y_0)$$

\no{and}

$$\rho_p(x_0,y_0) \leq \rho_p(x_n,x_0)+\rho_p(x_n,y_n)+\rho_p(y_n,y_0).$$

\no Therefore,

$$|\rho_p(x_n,y_n)-\rho_p(x_0,y_0)| \ \leq \ \rho_p(x_n,x_0)+\rho_p(y_n,y_0).$$

Now for any number $\epsilon > 0$, since $\lim\limits_{n}x_n =
x_0$ and $\lim\limits_{n}y_n = y_0$, there exist numbers
$N_1(\epsilon), N_1(\epsilon)\geq N_1$ and $N_2(\epsilon),
N_2(\epsilon)\geq N_2$ such that $\rho_p(x_n,x_0) \ \leq \
\frac{\epsilon}{2}$ if $n\geq N_1(\epsilon)$ and $\rho_p(y_n,y_0)
\ \leq \ \frac{\epsilon}{2}$ if $n\geq N_2(\epsilon)$. Whence, if
we choose $n\geq max\{N_1(\epsilon), N_2(\epsilon)\}$, then

$$|\rho_p(x_n,y_n)-\rho_p(x_0,y_0)| \ < \ \epsilon.\quad\quad \natural$$

Whether can a convergent sequence have more than one limiting
points? The following result answers this question.

\vskip 4mm

\no{\bf Theorem $1.4.8$} \ {\it If $\{x_n\}$ is a convergent
sequence in a multi-metric space
$\widetilde{M}=\bigcup\limits_{i=1}^m M_i$, then $\{x_n\}$ has
only one limit point.}

\vskip 3mm

{\it Proof} \ According to Theorem $1.4.6$, there exist integers
$N$ and $i, 1\leq i\leq m$ such that $x_n\in M_i$ if $n\geq N$.
Now if

$$\lim\limits_nx_n = x_1 \ {\rm and} \ \lim\limits_nx_n = x_2,$$

\no and $n\geq N$, by definition,

$$0\leq\rho_i(x_1,x_2)\leq\rho_i(x_n,x_1)+\rho_i(x_n,x_2).$$

\no Whence, we get that $\rho_i(x_1,x_2)=0$. Therefore, $x_1=x_2$.
\quad\quad $\natural$

\vskip 4mm

\no{\bf Theorem $1.4.9$} \ {\it Any convergent sequence in a
multi-metric space is a bounded points set.}

\vskip 3mm

{\it Proof} \ According to Theorem $1.4.8$, we obtain this result
immediately. \quad\quad $\natural$

A sequence $\{x_n\}$ in a multi-metric space
$\widetilde{M}=\bigcup\limits_{i=1}^m M_i$ is called a {\it Cauchy
sequence} if for any number $\epsilon >0$, there exist integers
$N(\epsilon)$ and $s, 1\leq s\leq m$ such that for any integers
$m,n\geq N(\epsilon)$, $\rho_s(x_m,x_n) \ < \epsilon$.

\vskip 4mm

\no{\bf Theorem $1.4.10$} \ {\it A Cauchy sequence $\{x_n\}$ in a
multi-metric space $\widetilde{M}=\bigcup\limits_{i=1}^m M_i$ is
convergent if and only if $|\{x_n\}\bigcap M_k|$ is finite or
infinite but $\{x_n\}\bigcap M_k$ is convergent in $(M_k,\rho_k)$
for $\forall k, 1\leq k\leq m$. }

\vskip 3mm

{\it Proof} \ The necessity of these conditions in this theorem is
known by Theorem $1.4.6$.

Now we prove the sufficiency. By definition, there exist integers
$s, 1\leq s\leq m$ and $N_1$ such that $x_n\in M_s$ if $n\geq
N_1$. Whence, if $|\{x_n\}\bigcap M_k|$ is infinite and
$\lim\limits_n\{x_n\}\bigcap M_k = x$, then there must be $k=s$.
Denote by $\{x_n\}\bigcap M_k = \{x_{k1},x_{k2},\cdots ,x_{kn},
\cdots\}$.

For any positive number $\epsilon >0$, there exists an integer
$N_2, N_2\geq N_1$ such that $\rho_k(x_m,x_n) \ < \
\frac{\epsilon}{2}$ and $\rho_k(x_{kn},x) \ < \
\frac{\epsilon}{2}$ if $m,n\geq N_2$. According to Theorem
$1.4.7$, we get that

$$\rho_k(x_n,x)\leq \rho_k(x_n,x_{kn})+\rho_k(x_{kn},x)\ < \ \epsilon$$

\no if $n\geq N_2$. Whence, $\lim\limits_nx_n = x. \quad\quad
\natural$

A multi-metric space $\widetilde{M}$ is said to be {\it completed}
if its every Cauchy sequence is convergent. For a completed
multi-metric space, we obtain two important results similar to
Theorems $1.4.3$ and $1.4.4$ in metric spaces.

\vskip 4mm

\no{\bf Theorem $1.4.11$} \ {\it Let
$\widetilde{M}=\bigcup\limits_{i=1}^m M_i$ be a completed
multi-metric space. For an $\epsilon$-disk sequence
$\{B(\epsilon_n,x_n)\}$, where $\epsilon_n > 0$ for $n=1,2,3,
\cdots$, if the following conditions hold:}

($i$) \ {\it $B(\epsilon_1,x_1)\supset B(\epsilon_2,x_2)\supset
B(\epsilon_3,x_3)\supset\cdots\supset
B(\epsilon_n,x_n)\supset\cdots$};

($ii$) \ {\it $\lim\limits_{n\to+\infty}\epsilon_n=0$},

\no {\it then $\bigcap\limits_{n=1}^{+\infty}B(\epsilon_n,x_n)$
only has one point. }

\vskip 3mm

{\it Proof} \ First, we prove that the sequence $\{x_n\}$ is a
Cauchy sequence in $\widetilde{M}$. By the condition $(i)$, we
know that if $m\geq n$, then $x_m\in B(\epsilon_m, x_m)\subset
B(\epsilon_n, x_n)$. Whence $\rho_i(x_m,x_n) \ < \epsilon_n$
provided $x_m,x_n\in M_i$ for $\forall i, 1\leq i\leq m$.

Now for any positive number $\epsilon$, since $\lim\limits_{n\to
+\infty}\epsilon_n=0$, there exists an integer $N(\epsilon)$ such
that if $n\geq N(\epsilon)$, then $\epsilon_n < \epsilon$.
Therefore, if $x_n\in M_l$, then
$\lim\limits_{m\rightarrow+\infty} x_m=x_n$. Thereby there exists
an integer $N$ such that if $m\geq N$, then $x_m\in M_l$ by
Theorem $1.4.6$. Choice integers $m,n\geq max\{N,N(\epsilon)\}$,
we know that

$$\rho_l(x_m,x_n) \ <  \ \epsilon_n \ < \ \epsilon.$$

\no So $\{x_n\}$ is a Cauchy sequence.

By the assumption that $\widetilde{M}$ is completed, we know that
the sequence $\{x_n\}$ is convergent to a point $x_0,
x_0\in\widetilde{M}$. By conditions of ($i$) and ($ii$), we get
that $\rho_l(x_0,x_n) \ < \epsilon_n$ if $m\to +\infty$. Whence,
$x_0\in\bigcap\limits_{n=1}^{+\infty}B(\epsilon_n,x_n)$.

Now if there is a point
$y\in\bigcap\limits_{n=1}^{+\infty}B(\epsilon_n,x_n)$, then there
must be $y\in M_l$. We get that

$$0\leq\rho_l(y,x_0) = \lim\limits_n\rho_l(y,x_n)\leq
\lim\limits_{n\to +\infty}\epsilon_n=0$$

\no by Theorem $1.4.7$. Therefore, $\rho_l(y,x_0)=0$. By the
definition of a metric function, we get that $y=x_0$. \quad\quad
$\natural$

Let $\widetilde{M}_1$ and $\widetilde{M}_2$ be two multi-metric
spaces and let $f: \widetilde{M}_1\rightarrow\widetilde{M}_2$ be a
mapping, $x_0\in\widetilde{M}_1, f(x_0)=y_0$. For $\forall\epsilon
> 0$, if there exists a number $\delta$ such that
$f(x)=y\in B(\epsilon ,y_0)\subset\widetilde{M}_2$ for $\forall
x\in B(\delta , x_0)$, i.e.,

$$f(B(\delta , x_0))\subset B(\epsilon ,y_0),$$

\no then we say that $f$ is {\it continuous at point} $x_0$. A
mapping $f:\widetilde{M}_1\rightarrow\widetilde{M}_2$ is called a
{\it continuous mapping} from $\widetilde{M}_1$ to
$\widetilde{M}_2$ if $f$ is continuous at every point of
$\widetilde{M}_1$.

For a continuous mapping $f$ from $\widetilde{M}_1$ to
$\widetilde{M}_2$ and a convergent sequence $\{x_n\}$ in
$\widetilde{M}_1$, $\lim\limits_nx_n=x_0$, we can prove that

$$\lim\limits_nf(x_n)=f(x_0).$$

For a multi-metric space $\widetilde{M}=\bigcup\limits_{i=1}^m
M_i$ and a mapping $T: \widetilde{M}\rightarrow\widetilde{M}$, if
there is a point $x^*\in\widetilde{M}$ such that $Tx^* = x^*$,
then $x^*$ is called a {\it fixed point} of $T$. Denote the number
of fixed points of a mapping $T$ in $\widetilde{M}$ by $^{\#}\Phi
(T)$. A mapping $T$ is called a {\it contraction} on a
multi-metric space $\widetilde{M}$ if there are a constant $\alpha
, 0 < \alpha < 1$ and integers $i, j, 1\leq i, j\leq m$ such that
for $\forall x,y\in M_i$, $Tx, Ty\in M_j$ and

$$\rho_j(Tx,Ty)\leq\alpha\rho_i(x,y).$$

\vskip 4mm

\no{\bf Theorem $1.4.12$} \ {\it Let
$\widetilde{M}=\bigcup\limits_{i=1}^m M_i$ be a completed
multi-metric space and let $T$ be a contraction on
$\widetilde{M}$. Then}

$$1\leq ^{\#}\Phi (T)\leq m.$$

\vskip 3mm

{\it Proof} \ Choose arbitrary points $x_0, y_0\in M_1$ and define
recursively

$$ x_{n+1}=Tx_n, \ \ y_{n+1}=Tx_n$$

\no for $n=1,2,3,\cdots$. By definition, we know that for any
integer $n, n\geq 1$, there exists an integer $i, 1\leq i\leq m$
such that $x_n, y_n\in M_i$. Whence, we inductively get that

$$0\leq\rho_i(x_n,y_n)\leq\alpha^n\rho_1(x_0,y_0).$$

Notice that $0 < \alpha < 1$, we know that $\lim\limits_{n\to
+\infty}\alpha^n=0$. Thereby there exists an integer $i_0$ such
that

$$\rho_{i_0}(\lim\limits_nx_n,\lim\limits_ny_n)=0.$$

\no Therefore, there exists an integer $N_1$ such that $x_n,
y_n\in M_{i_0}$ if $n\geq N_1$. Now if $n\geq N_1$, we get that

\begin{eqnarray*}
\rho_{i_0}(x_{n+1},x_n)&=& \rho_{i_0}(Tx_n,Tx_{n-1})\\
&\leq& \alpha
\rho_{i_0}(x_n,x_{n-1})=\alpha\rho_{i_0}(Tx_{n-1},Tx_{n-2})\\
&\leq& \alpha^2
\rho_{i_0}(x_{n-1},x_{n-2})\leq\cdots\leq\alpha^{n-N_1}\rho_{i_0}(x_{N_1+1},x_{N_1}).
\end{eqnarray*}

\no and generally, for $m\geq n\geq N_1$,

\begin{eqnarray*}
\rho_{i_0}(x_m,x_n) &\leq&
\rho_{i_0}(x_n,x_{n+1})+\rho_{i_0}(x_{n+1},x_{n+2})+\cdots
+\rho_{i_0}(x_{n-1},x_n)\\
&\leq& (\alpha^{m-1}+\alpha^{m-2}+\cdots
+\alpha^{n})\rho_{i_0}(x_{N_1+1},x_{N_1})\\
&\leq& \frac{\alpha^n}{1-\alpha}\rho_{i_0}(x_{N_1+1},x_{N_1})\to 0
(m,n\to +\infty).
\end{eqnarray*}

\no Therefore, $\{x_n\}$ is a Cauchy sequence in $\widetilde{M}$.
Similarly, we can also prove $\{y_n\}$ is a Cauchy sequence.

Because $\widetilde{M}$ is a completed multi-metric space, we know
that

$$\lim\limits_nx_n = \lim\limits_ny_n = z^*.$$

Now we prove $z^*$ is a fixed point of $T$ in $\widetilde{M}$. In
fact, by $\rho_{i_0}(\lim\limits_nx_n,\lim\limits_ny_n)=0$, there
exists an integer $N$ such that

$$x_n,y_n, Tx_n,Ty_n\in M_{i_0}$$

\no if $n\geq N+1$. Whence, we know that

\begin{eqnarray*}
0\leq\rho_{i_0}(z^*,Tz^*) &\leq& \rho_{i_0}(z^*,x_n)+
\rho_{i_0}(y_n,Tz^*)+\rho_{i_0}(x_n,y_n)\\
&\leq& \rho_{i_0}(z^*,x_n)+\alpha
\rho_{i_0}(y_{n-1},z^*)+\rho_{i_0}(x_n,y_n).
\end{eqnarray*}

\no Notice that

$$\lim\limits_{n\to +\infty}\rho_{i_0}(z^*,x_n)=\lim\limits_{n\to
+\infty}\rho_{i_0}(y_{n-1},z^*)=\lim\limits_{n\to
+\infty}\rho_{i_0}(x_n,y_n)=0.$$

\no We get $\rho_{i_0}(z^*,Tz^*)=0$, i.e., $Tz^*=z^*$.

For other chosen points $u_0,v_0\in M_1$, we can also define
recursively

$$u_{n+1} = Tu_n, \ \ v_{n+1} = Tv_n$$

\no and get a limiting point $\lim\limits_nu_n = \lim\limits_nv_n
= u^*\in M_{i_0}, Tu^*\in M_{i_0}$. Since

$$\rho_{i_0}(z^*,u^*)=\rho_{i_0}(Tz^*,Tu^*)\leq\alpha\rho_{i_0}(z^*,u^*)$$

\no and $0 < \alpha < 1,$ there must be $z^*=u^*$.

Similarly consider the points in $M_i, 2\leq i\leq m$, we get that

$$1\leq ^{\#}\Phi (T)\leq m. \quad\quad\quad \natural$$

As a consequence, we get the {\it Banach theorem} in metric
spaces.

\vskip 4mm

\no{\bf Corollary $1.4.2$}(Banach) \ {\it Let $M$ be a metric
space and let $T$ be a contraction on $M$. Then $T$ has just one
fixed point.}

\vskip 8mm

\no{\bf \S $1.5$ \ Remarks and Open Problems}

\vskip 6mm

\no The central idea of Smarandache multi-spaces is to combine
different fields (spaces, systems, objects, $\cdots$) into a
unifying field and find its behaviors. Which is entirely new, also
an application of combinatorial approaches to classical
mathematics but more important than combinatorics itself. This
idea arouses us to think why an assertion is true or not in
classical mathematics. Then combine an assertion with its
non-assertion and enlarge the filed of truths. A famous fable says
that {\it each theorem in mathematics is an absolute truth}. But
we do not think so. Our thinking is that {\it each theorem in
mathematics is just a relative truth}. Thereby we can establish
new theorems and present new problems boundless in mathematics.
Results obtained in Section $1.3$ and $1.4$ are applications of
this idea to these groups, rings, vector spaces or metric spaces.
Certainly, more and more multi-spaces and their good behaviors can
be found under this thinking. Here we present some remarks and
open problems for multi-spaces.

\vskip 3mm

\no{\bf $1.5.1.$ Algebraic Multi-Spaces} \ The algebraic
multi-spaces are discrete representations for phenomena in the
natural world. They maybe completed or not in cases. For a
completed algebraic multi-space, it is a reflection of an
equilibrium phenomenon. Otherwise, a reflection of a
non-equilibrium phenomenon. Whence, more consideration should be
done for algebraic multi-spaces, especially, by an analogous
thinking as in classical algebra.

\vskip 3mm

\no{\bf Problem $1.5.1$} \ {\it Establish a decomposition theory
for multi-groups.}

\vskip 2mm

In group theory, we know the following decomposition
result([$107$][$82$]) for groups.

\vskip 2mm

{\it Let $G$ be a finite $\Omega$-group. Then $G$ can be uniquely
decomposed as a direct product of finite non-decomposition
$\Omega$-subgroups.}

\vskip 2mm

{\it Each finite abelian group is a direct product of its Sylow
$p$-subgroups.}\vskip 2mm

Then Problem $1.5.1$ can be restated as follows.

\vskip 2mm

{\it Whether can we establish a decomposition theory for
multi-groups similar to the above two results in group theory,
especially, for finite multi-groups?}

\vskip 3mm

\no{\bf Problem $1.5.2$} \ {\it Define the conception of simple
multi-groups. For finite multi-groups, whether can we find all
simple multi-groups?}\vskip 2mm

For finite groups, we know that there are four simple group
classes ([$108$]):

\vskip 2mm

{\bf Class $1$}: the cyclic groups of prime order; \vskip 2mm

{\bf Class $2$}: the alternating groups $A_n, n\geq 5$;

\vskip 2mm

{\bf Class $3$}: the 16 groups of Lie types;

\vskip 2mm

{\bf Class $4$}: the 26 sporadic simple groups.

\vskip 3mm

\no{\bf Problem $1.5.3$} \ {\it Determine the structure properties
of multi-groups generated by finite elements.} \vskip 2mm

For a subset $A$ of a multi-group $\widetilde{G}$, define its
spanning set by

$$\left<A\right>=\{a\circ b | a,b\in A \ {\rm and} \ \circ\in
O(\widetilde{G})\}.$$

\no If there exists a subset $A\subset\widetilde{G}$ such that
$\widetilde{G}=\left<A\right>$, then call $\widetilde{G}$ is
generated by $A$. Call $\widetilde{G}$ is {\it finitely generated}
if there exist a finite set $A$ such that
$\widetilde{G}=\left<A\right>$. Then Problem $5.3$ can be restated
by \vskip 2mm

{\it Can we establish a finite generated multi-group theory
similar to the finite generated group theory?}

\vskip 3mm

\no{\bf Problem $1.5.4$} \ {\it Determine the structure of a
Noether multi-ring. }

\vskip 2mm

Let $R$ be a ring. Call $R$ a {\it Noether ring} if its every
ideal chain only has finite terms. Similarly, for a multi-ring
$\widetilde{R}$, if its every multi-ideal chain only has finite
terms, it is called a {\it Noether multi-ring}. {\it Whether can
we find its structures similar to Corollary $1.3.5$ and Theorem
$1.3.12$?}

\vskip 3mm

\no{\bf Problem $1.5.5$} \ {\it Similar to ring theory, define a
Jacobson or Brown-McCoy radical for multi-rings and determine
their contribution to multi-rings.}\vskip 2mm

Notice that Theorem $1.3.14$ has told us there is a similar linear
theory for multi-vector spaces, but the situation is more complex.

\vskip 3mm

\no{\bf Problem $1.5.6$} \ {\it Similar to linear spaces, define
linear transformations on multi-vector spaces. Can we establish a
matrix theory for these linear transformations?}

\vskip 3mm

\no{\bf Problem $1.5.7$} \ {\it Whether a multi-vector space must
be a linear space?}

\vskip 2mm

\no{\bf Conjecture $1.5.1$} \ {\it There are non-linear
multi-vector spaces in multi-vector spaces.}\vskip 2mm

Based on Conjecture $1.5.1$, there is a fundamental problem for
multi-vector spaces.

\no{\bf Problem $1.5.8$} \ {\it Can we apply multi-vector spaces
to non-linear spaces?}

\vskip 3mm

\no{\bf $1.5.2.$ Multi-Metric Spaces} \ On a tradition notion,
only one metric maybe considered in a space to ensure the same on
all the time and on all the situation. Essentially, this notion is
based on an assumption that all spaces are homogeneous. In fact,
it is not true in general.

Multi-metric spaces can be used to simplify or beautify
geometrical figures and algebraic equations. For an explanation,
an example is shown in Fig.$1.3$, in where the left elliptic curve
is transformed to the right circle by changing the metric along
$x,y$-axes and an elliptic equation

$$\frac{x^2}{a^2}+\frac{y^2}{b^2}=1$$

\no to equation

$$x^2+y^2=r^2$$

\no of a circle of radius $r$.

\vskip 2mm

\includegraphics[bb=75 5 400 170]{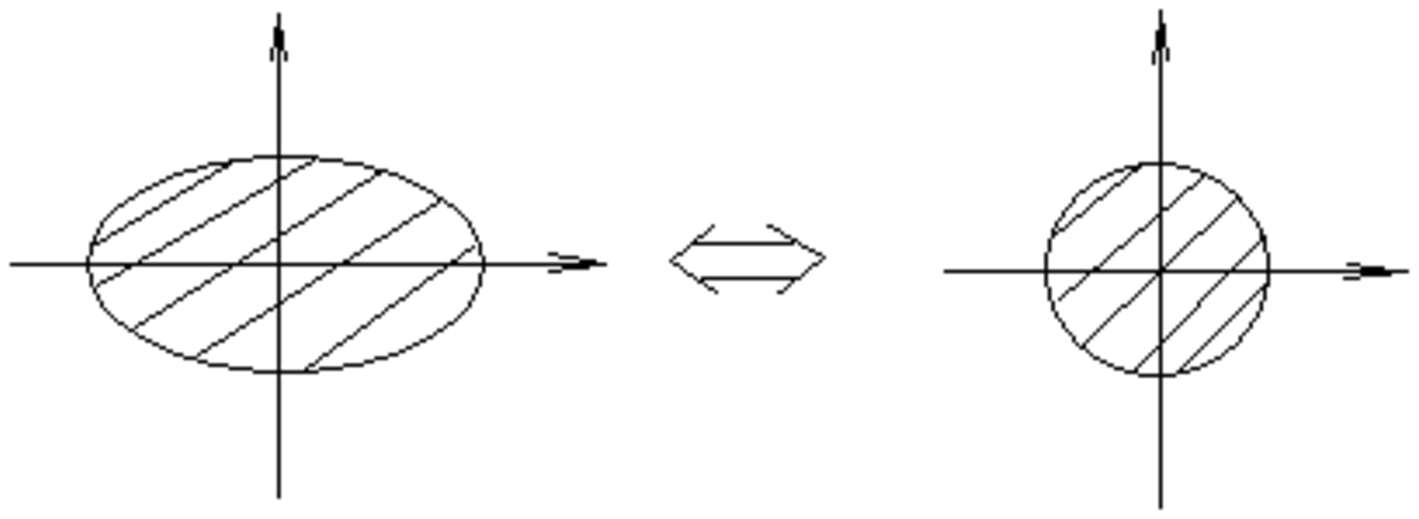}

\vskip 2mm

\c{\bf Fig.$1.3$}\vskip 2mm

Generally, in a multi-metric space, we can simplify a polynomial
similar to the approach used in the projective geometry. {\it
Whether this approach can be contributed to mathematics with
metrics?}

\no{\bf Problem $1.5.9$} \ {\it Choose suitable metrics to
simplify the equations of surfaces or curves in ${\bf R}^3$.}

\vskip 3mm

\no{\bf Problem $1.5.10$} \ {\it Choose suitable metrics to
simplify the knot problem. Whether can it be used for classifying
$3$-dimensional manifolds?}

\vskip 3mm

\no{\bf Problem $1.5.11$} \ {\it Construct multi-metric spaces or
non-linear spaces by Banach spaces. Simplify equations or problems
to linear problems.}

\vskip 4mm

\no{\bf $1.5.3.$ Multi-Operation Systems} \ By a complete
Smarandache multi-space $\widetilde{A}$ with an operation set
$O(\widetilde{A})$, we can get a {\it multi-operation system}
$\widetilde{A}$. For example, if $\widetilde{A}$ is a multi-field
$\widetilde{F}=\bigcup\limits_{i=1}^nF_i$ with an operation set
$O(\widetilde{F})=\{(+_i,\times_i) | \ 1\leq i\leq n\}$, then
$(\widetilde{F}; +_1,+_2,\cdots ,+_n)$, $(\widetilde{F};
\times_1,\times_2,\cdots ,\times_n)$ and
$(\widetilde{F};(+_1,\times_1),$ $(+_2,\times_2),\cdots
,(+_n,\times_n))$ are multi-operation systems. On this view, the
classical operation system $(R \ ;+)$ and $(R \ ;\times )$ are
only {\it sole operation systems}. For a multi-operation system
$\widetilde{A}$, we can define these conceptions of equality and
inequality, $\cdots$, etc.. For example, in the multi-operation
system $(\widetilde{F}; +_1,+_2,\cdots ,+_n)$, we define the
equalities $=_1,=_2,\cdots , =_n$ such as those in sole operation
systems $(\widetilde{F};+_1),(\widetilde{F};+_2),\cdots ,$
$(\widetilde{F};+_n)$, for example, $2=_1 2, 1.4=_2 1.4, \cdots ,
\sqrt{3}=_n \sqrt{3}$ which is the same as the usual meaning and
similarly, for the conceptions $\geq_1,\geq_2,\cdots ,\geq_n$ and
$\leq_1,\leq_2,\cdots ,\leq_n$.

In a classical operation system $(R \ ;+)$, the equation system

\begin{eqnarray*}
x+2+4+6 &=& 15\\
x+1+3+6 &=& 12\\
x+1+4+7 &=& 13
\end{eqnarray*}

\no can not has a solution. But in the multi-operation system
$(\widetilde{F}; +_1,+_2,\cdots ,+_n)$, the equation system

\begin{eqnarray*}
x+_12+_14+_16 &=_1& 15\\
x+_21+_23+_26 &=_2& 12\\
x+_31+_34+_37 &=_3& 13
\end{eqnarray*}

\no may have a solution $x$ if

\begin{eqnarray*}
15+_1(-1)+_1(-4)+_1(-16) &=&
12+_2(-1)+_2(-3)+_2(-6)\\
&=& 13+_3(-1)+_3(-4)+_3(-7).
\end{eqnarray*}

\no in $(\widetilde{F}; +_1,+_2,\cdots ,+_n)$. Whence, an element
maybe have different disguises in a multi-operation system.

For the multi-operation systems, a number of open problems needs
to research further.

\vskip 4mm

\no{\bf Problem $1.5.12$} \ {\it Find necessary and sufficient
conditions for a multi-operation system with more than $3$
operations to be the rational number field $Q$, the real number
field $R$ or the complex number field $C$.}

\vskip 3mm

For a multi-operation system $(N \ ;
(+_1,\times_1),(+_2,\times_2),\cdots ,(+_n,\times_n))$ and
integers $a,b,c\in N$, if $a=b\times_ic$ for an integer $i, 1\leq
i\leq n$, then $b$ and $c$ are called {\it factors} of $a$. An
integer $p$ is called a {\it prime} if there exist integers
$n_1,n_2$ and $i, 1\leq i\leq n$ such that $p=n_1\times_i n_2$,
then $p=n_1$ or $p=n_2$. Two problems for primes of a
multi-operation system $(N \ ; (+_1,\times_1),$
$(+_2,\times_2),\cdots ,(+_n,\times_n))$ are presented in the
following.

\vskip 4mm

\no{\bf Problem $1.5.13$} \ {\it For a positive real number $x$,
denote by $\pi_m(x)$ the number of primes $\leq x$ in $(N \ ;
(+_1,\times_1),(+_2,\times_2),\cdots ,(+_n,\times_n))$. Determine
or estimate $\pi_m(x)$.}

Notice that for the positive integer system, by a well-known
theorem, i.e., {\it Gauss prime theorem}, we have known
that([$15$])

$$\pi (x)\sim\frac{x}{{\rm log}x}.$$

\vskip 3mm

\no{\bf Problem $1.5.14$} \ {\it Find the additive number
properties for $(N \ ; (+_1,\times_1),(+_2,\times_2),\cdots ,$
$(+_n,\times_n))$, for example, we have weakly forms for
Goldbach's conjecture and Fermat's problem {\rm ([$34$])} as
follows.

\vskip 3mm

\no{\bf Conjecture $1.5.2$} \ For any even integer $n, n\geq 4$,
there exist odd primes $p_1,p_2$ and an integer $i, 1\leq i\leq n$
such that $n=p_1+_ip_2$.

\vskip 3mm

\no{\bf Conjecture $1.5.3$} \ For any positive integer $q$, the
Diophantine equation $x^q+y^q=z^q$ has non-trivial integer
solutions $(x,y,z)$ at least for an operation ¡°$+_i$¡± with
$1\leq i\leq n$. }

\vskip 3mm

A {\it Smarandache $n$-structure on a set $S$} means a weak
structure $\{w(0)\}$ on $S$ such that there exists a chain of
proper subsets $P(n-1)\subset P(n-2)\subset\cdots\subset
P(1)\subset S$ whose corresponding structures verify the inverse
chain
$\{w(n-1)\}\supset\{w(n-2)\}\supset\cdots\supset\{w(1)\}\supset\{w(0)\}$,
i.e., structures satisfying more axioms.

\vskip 4mm

\no{\bf Problem $1.5.15$} \ {\it For Smarandache multi-structures,
solves these Problems $1.5.1-1.5.8$.}\vskip 3mm

\no{\bf $1.5.4.$ Multi-Manifolds} \ Manifolds are important
objects in topology, Riemann geometry and modern mechanics. It can
be seen as a local generalization of Euclid spaces. By the
Smarandache's notion, we can also define multi-manifolds. To
determine their behaviors or structure properties will useful for
modern mathematics.

In an Euclid space ${\bf R}^n$, an {\it $n$-ball of radius $r$} is
defined by

$$B^n(r) = \{(x_1,x_2,\cdots ,x_n) | x_1^2+x_2^2+\cdots +x_n^2\leq r \}.$$

Now we choose $m$ $n$-balls $B_1^n(r_1), B_2^n(r_2), \cdots
,B_m^n(r_m)$, where for any integers $i,j, 1\leq i,j\leq m$,
$B_i^n(r_i)\bigcap B_j^n(r_j)=\empty$ or not and $r_i=r_j$ or not.
An {\it $n$-multi-ball} is a union

$$\widetilde{B} = \bigcup\limits_{k=1}^mB_k^n(r_k).$$

\no Then an {\it $n$-multi-manifold} is a Hausdorff space with
each point in this space has a neighborhood homeomorphic to an
$n$-multi-ball.

\vskip 3mm

\no{\bf Problem $1.5.16$} \ {\it For an integer $n, n\geq 2$,
classifies $n$-multi-manifolds. Especially, classifies
$2$-multi-manifolds.}

\vskip 2mm

For closed $2$-manifolds, i.e., locally orientable surfaces, we
have known a classification theorem for them.

\vskip 3mm

\no{\bf Problem $1.5.17$} \ {\it If we replace the word
¡°homeomorphic¡± by ¡°points equivalent¡± or ¡°isomorphic¡±, what
can we obtain for $n$-multi-manifolds? Can we classify them? }

\vskip 2mm

Similarly, we can also define differential multi-manifolds and
consider their contributions to modern differential geometry,
Riemann geometry or modern mechanics, $\cdots$, etc..

\end{document}